\documentclass{amsart}
\usepackage{amssymb, latexsym}
\def\be#1{\begin{equation} \label{#1}}
\def\bi{\begin{itemize}}
\def\bs{\begin{split}}
\def\es{\end{split}}
\def\ba{\begin{align}}
\def\bas{\begin{align*}}
\def\ea{\end{align}}
\def\eas{\end{align*}}

\def\sgn{{\hbox{sgn}}}

\def\R{{{\mathbb R}}}
\def\P{{{\bf P}}}
\def\Z{{{\mathbb Z}}}
\def\T{{{\mathbb T}}}

\def\eps{\varepsilon}
\newcommand{\half}{\frac{1}{2}}
\def \endprf{\hfill  {\vrule height6pt width6pt depth0pt}\medskip}
\def\emph#1{{\it #1}}
\def\textbf#1{{\bf #1}}
\newtheorem{theorem}{Theorem}
\newtheorem{proposition}{Proposition}
\newtheorem{lemma}{Lemma}
\newtheorem{corollary}{Corollary}
\theoremstyle{definition}
\newtheorem{definition}{Definition}
\theoremstyle{remark}
\newtheorem{remark}{Remark}
\numberwithin{equation}{section}
\numberwithin{lemma}{section}
\numberwithin{remark}{section}

\begin{document}

\title[Multilinear estimates and applications]{Multilinear estimates for periodic KdV equations, and applications}
\author{J.~Colliander}
\thanks{J.E.C. is supported in part by N.S.F. grant DMS 0100595.}
\address{University of Toronto }
\author{M.~Keel}
\thanks{M.K. is supported in part by N.S.F. Grant DMS 9801558}
\address{University of Minnesota}
\author{G.~Staffilani}
\thanks{G.S. is supported in part by N.S.F. Grant 9800879 and by grants 
from Hewlett and Packard and the Sloan Foundation.}
\address{Brown University}
\author{H.~Takaoka}
\thanks{H.T. is supported in part by J.S.P.S. Grant No. 13740087.}
\address{Hokkaido University}
\author{T.~Tao}
\thanks{T.T. is a Clay Prize Fellow and is supported in part by grants
from the Packard Foundation.}
\address{University of California, Los Angeles}

%\vspace{-0.3in}

\begin{abstract}
We prove an endpoint multilinear estimate for the $X^{s,b}$ spaces associated to the periodic Airy equation.  As a consequence we obtain sharp local well-posedness results for periodic generalized KdV equations, as well as some global well-posedness results below the energy norm.
\end{abstract}

\date{8 February 2002}

\maketitle

\tableofcontents

\section{Introduction}\label{introduction-sec}

This paper studies the 
Cauchy problem for periodic generalized KdV equations of the form
\begin{equation}
\label{pkdv}
\left\{
   \begin{matrix}
    \partial_t u + \frac{1}{4\pi^2} \partial_x^3 u +  F(u)_x =0,& u: \T \times 
[0,T] \longmapsto \R ,
         \\
     u( x, 0) = \phi(x), & x \in \T 
   \end{matrix}
\right.
\end{equation}
where $F$ is a polynomial of degree $k+1$, the initial data $u_0$ is in a  Sobolev space $H^s (\T) $, and $\T = \R / \Z $ is the torus. 
 The factor $\frac{1}{4\pi^2}$ is convenient in order to make the dispersion relation $\tau = \xi^3$, but it is inessential and we recommend that the reader ignore all powers of $2\pi$ which appear in the sequel.  We can assume that $F$ has no constant or linear term since these can be removed by a Gallilean transformation.

The main result established in this paper is a sharp multilinear estimate
which allows us to show the initial value problem \eqref{pkdv} is locally
well-posed in $H^s ( \T )$ for $s \geq \half$. We expect that the estimate
(contained in Theorem \ref{main} below) will have other applications in the
study of the behavior of solutions of KdV-like equations. 

If $u(x,t)$ is a function on the cylinder $\T \times \R$ and $s,b \in \R$, we define the $X^{s,b} = X^{s,b}_{\tau = \xi^3}(\T \times \R)$ norm by
$$
\| u \|_{X^{s,b}} := \| \hat u(\xi, \tau) \langle \xi \rangle^s \langle \tau - \xi^3\rangle^b \|_{L^2_{\tau,\xi}}
$$
where the space-time Fourier transform $\hat u(\xi,\tau)$ is given for $\xi \in \Z$, $\tau \in \R$ by
$$ \hat u(\xi, \tau) := \int_{\T \times \R} e^{-2 \pi i (x \xi + t\tau)} u(x,t)\ dx dt
$$
and $\langle x \rangle := 1 + |x|$.  We shall often abbreviate $\| u\|_{X^{s,b}}$ as $\| u \|_{s,b}$.  These norms were introduced for the KdV equation in \cite{borg:xsb} (with similar spaces for the wave equation appearing in \cite{beals:xsb}, \cite{kl-mac:xsb}).

The low-regularity study of the equation \eqref{pkdv} on $\T$ has been based around iteration in the spaces $X^{s,\frac{1}{2}}$ (see \cite{borg:xsb}, \cite{kpv:kdv}, \cite{staff:pkdv}).  This space barely fails to control the $L^\infty_t H^s_x$ norm.  To rectify this we define the slightly stronger norms $Y^s$ by
\be{y-def}
\| u \|_{Y^s} := \| u \|_{s,\frac{1}{2}} + \| \langle \xi \rangle^s \hat u \|_{L^2_\xi L^1_\tau}.
\end{equation}

We shall also need the companion spaces $Z^s$ defined by
\be{z-def}
\| u \|_{Z^s} := \| u \|_{s,-\half} + \left\| 
\frac{\langle \xi \rangle^s \hat u}{\langle \tau - \xi^3  \rangle} \right\|_{L^2_\xi L^1_\tau}.
\end{equation}

\subsection{Main Results}

The main technical result of this paper is the following multilinear estimate.
\begin{theorem}\label{main}
For any $s \geq \half$, we have\footnote{Here and in the sequel, $A \lesssim B$ denotes $A \leq C B$, where $C$ is a constant possibly depending on $s$, $k$.}
\be{main-est}
\| \prod_{i=1}^k u_i \|_{s-1,\half} \lesssim \prod_{i=1}^k \| u_i \|_{Y^s}.
\end{equation}
\end{theorem}

This improves on the results in \cite{staff:pkdv}, where this estimate was proven for $s \geq 1$.

By combining this estimate with an estimate of Kenig, Ponce, and Vega \cite{kpv:kdv} we shall obtain
\begin{proposition}\label{main-cor}
For any $s \geq \half$ and functions $u_1, \ldots, u_{k+1}$ in $Y^s$ we have
\begin{equation}
\label{star} 
\left\| \P(\P(\prod_{i=1}^k u_i) \partial_x u_{k+1}) \right\|_{Z^s} \lesssim \prod_{i=1}^{k+1} \| u_i \|_{Y^s},
\end{equation}
where $\P$ denotes the orthogonal projection onto mean zero functions
$$ \P(u)(x) := u(x) - \int_\T u.$$
Moreover, the estimate \eqref{star} fails for $s< \frac{1}{2}.$
\end{proposition}

In \cite{kpv:kdv} the restriction $s \geq \half$ was shown to be sharp in the $k=2$ case.  The counter-example given there can be easily modified to show the necessity of the condition $s \geq \half$ for $k \geq 4$, and also for the $k=3$ case if one allows the mean of the $u_i$ to be non-zero.  In the $k=3$ case with mean zero the counter-example is a little trickier, and is discussed in Section \ref{counter-sec}. 

The sharp estimate \eqref{star} of Proposition \ref{main-cor} will be useful in studying the dynamical behavior of solutions of polynomial generalizations of the KdV equation. In this paper, we use Proposition \ref{main-cor} to obtain some local and global well-posedness results for the Cauchy problem \eqref{pkdv}.  In \cite{staff:pkdv} (see also \cite{borg:xsb}, \cite{kpv:contract}), such equations were shown to be globally well-posed for $H^1$ data.

\begin{theorem}\label{periodic-kdv}  If $F(u)$ is a polynomial, then the Cauchy problem for the periodic generalized KdV equation \eqref{pkdv}
is locally well-posed in $H^s(\T)$ for all $s \geq \half$, if the $H^s$ norm of the data is sufficiently small. 
\end{theorem}

For the large data case, see Remarks \ref{largedata} and \ref{largedatatwo} below. When $F$ is quadratic we have the KdV equation, which is locally well-posed all the way down to $s = -\half$ \cite{kpv:kdv} (see also \cite{borg:xsb}, \cite{borg:measures}).  When $F$ is cubic we have the modified KdV equation, for which local well-posedness was already obtained for $s \geq \half$ in \cite{kpv:kdv} (see also \cite{ckstt:2}), and this range is sharp for the purposes of uniformly continuous dependence of the solution on the data, see \cite{kpv:kdv}, \cite{kpv:counter}.  In the quartic case $F(u) = u^4$ we shall show in Section \ref{counter-sec} that one has analytic ill-posedness in $H^s$ for any $s < \half$; this example can be extended to general polynomials $F$ of degree at least 3. 

On the real line with $F(u) = u^{k+1}$, the equation \eqref{pkdv} is known to be locally well-posed down to the scaling exponent $s \geq \half - \frac{2}{k}$ for $k \geq 4$ \cite{kpv:contract} (see also \cite{birnir}; earlier results are in \cite{gt}).  This was recently extended to $k=3$ (except at the endpoint $s = \half - \frac{2}{k}$) in \cite{axel}.  Thus there is a loss of $\frac{2}{k}$ derivatives when moving to the periodic setting when $k \geq 3$ in contrast to the $\frac{1}{4}$ loss one has in the $k=1,2$ cases. These observations resolve a problem\footnote{See the problem posed after Theorem 5.3 in Kenig's lecture notes at http://www.msri.org/publications/ln/msri/1997/ha/kenig/5/banner/03.html.} posed by Carlos Kenig.

There are some other consequences of Theorem \ref{main}.  It allows us to complete the proof (in \cite{ckstt:2}) of global well-posedness of the KdV and mKdV equations down to $s \geq -\half$ and $s \geq \half$.  In particular in the symplectic space $\dot H^{-\half}$ we see that the KdV flow is a smooth symplectic flow (see \cite{borg:book}, \cite{KuksinSqueeze} for further discussion).  Also, we can obtain some periodic global well-posedness results for generalized KdV equations below $H^1$.  Specifically, we have

\begin{theorem}\label{gwp}  The Cauchy problem for the periodic generalized KdV equation
\begin{equation}
\label{gkdv}
\left\{
\begin{matrix}
u_t + \frac{1}{4\pi^2} u_{xxx} + u^3 u_x = 0 \\
u(x,0) = u_0(x), ~ x \in \T
\end{matrix}
\right.
\end{equation}
is locally well-posed for large $H^s(\T)$ data for $s \geq \half$, and globally well-posed for large real $H^s(\T)$ data for $s > \frac{5}{6}$. Moreover, the initial value problem \eqref{gkdv} is analytically ill-posed in $H^s ( \T )$ for $s < \frac{1}{2}$.
\end{theorem}

We prove this in Sections \ref{period-sec}-\ref{global-sec}, by using the ``$I$-method'' in \cite{keel:mkg}, \cite{ckstt:1}, \cite{ckstt:2a} (see also \cite{keel:wavemap}, \cite{ckstt:2}); this is a method for obtaining global well-posedness below the energy norm by constructing ``almost conserved'' analogues of the Hamiltonian for rough solutions.  We remark that the truncation method of Bourgain (\cite{borg:book}, \cite{fonesca}, \cite{cst}, $\ldots$) to obtain global well-posedness below the energy norm does not apply here because the equation has no smoothing properties in the periodic context.

\subsection{Remarks and possible extensions}

As in many other results on global well posedness below the energy norm, our condition $s > \frac{5}{6}$ for global well-posedness falls quite short of the local condition $s \geq \half$.  This is mainly due to our poor control on the fluctuation of our modified Hamiltonian.  From our experience with the KdV and mKdV equations in \cite{ckstt:2} however we believe it is reasonable to hope that global well-posedness should hold for all $s \geq \half$.

One should be able to obtain an analogue of Theorem \ref{gwp} for defocusing generalized KdV equations (i.e. with the non-linearity $u^3 u_x$ replaced by $\alpha u^k u_x$ for $k$ even and $\alpha \leq 0$).  Morally speaking, our method gives this whenever $s > \frac{13}{14} - \frac{2}{7k}$ (see the footnote in Section \ref{global-sec}), however there are some technical difficulties when $k>4$ because the conserved Hamiltonian\footnote{Admittedly one also has conservation of $L^2$ norm, but this turns out to not be so useful because $L^2$ is super-critical when $k > 4$ and therefore does not scale favorably.}
$$\int \frac{1}{8 \pi^2} u_x^2 - \frac{\alpha}{(k+1)(k+2)} u^{k+2}\ dx$$
does not quite control the $H^1$ norm of $u$ due to low frequency issues (and the fact that $L^2$ is now supercritical).  In principle this could be avoided by the techniques in \cite{keel:mkg}, but we do not present this here, as in any event these results are almost certainly not sharp and should be significantly improvable by adding correction terms to the modified Hamiltonian (see \cite{ckstt:2}, \cite{ckstt:3}, \cite{ckstt:4}).  We remark that when $k=1$ or $k=2$ we can obtain global well-posedness for the same range of exponents as the local theory ($s \geq -\half$ and $s \geq \half$ respectively); see \cite{ckstt:2}.

We expect that Theorem \ref{gwp} may be generalized to include the case when $F$ is analytic or smooth by adapting the arguments in \cite{borg:gkdv}.

One should also be able to obtain similar global well-posedness results for the line $\R$ (with better exponents than the periodic case).  In fact, the arguments should be more elementary, requiring no number theory and relying instead on such estimates as the Kato smoothing estimate and the sharp maximal function estimate (cf. \cite{ckstt:2}, \cite{fonesca}, \cite{kpv:contract}).

\subsection{Outline}
Section 2 shows the ill-posedness claim of Theorem \ref{gwp} and, under suitable modifications, the failure of \eqref{star} for $s< \frac{1}{2}.$ Section 3 records linear estimates between spaces associated to the Airy equation. Section 4 reduces Theorem \ref{main} to a multiplier bound which we establish in the non-endpoint case, $s > \frac{1}{2}$, in Section 5. Section 6 recalls some elementary number theory which we use in Sections 7 and 8 to complete the proof of Theorem \ref{main} at the endpoint $s = \frac{1}{2}.$ Section 9 establishes Proposition \ref{main-cor}. Theorem \ref{periodic-kdv} is proven in Section 10. Section 11 rescales various estimates to the setting of large periods\footnote{Recasting the estimates in this form may be relevant in studying zero dispersive and semiclassical limit problems.}. Section 12 contains a general interpolation result revealing a certain flexibility in proofs of local well-posedness. The global result of Theorem \ref{gwp} is proven in Sections 13 and 14.

The authors thank Jorge Silva for a correction.

\section{A counter-example}\label{counter-sec}

In this section we give an example which shows why the condition $s \geq \half$ is necessary.  We shall discuss only the most difficult case, namely the $k=3$ case \eqref{gkdv} when $u$ has mean zero.  The other cases can be treated either by modifying the example given here or the one in \cite{kpv:kdv}.

Let $s < \half$, and let $N \gg 1$ be a large integer.  Let $N_0, N_1, N_2, N_3, N_4$ be integers with distinct magnitudes such that
$$ |N_0| \sim 1; \quad |N_1|, |N_2|, |N_3|, |N_4| \sim N$$
and
\be{whip}
N_0 + N_1 + N_2 + N_3 + N_4 = 0; \quad N_0^3 + N_1^3 + N_2^3 + N_3^3 + N_4^3 = O(1).
\end{equation}
For instance, we could choose
$$ N_0 := 6; \quad N_1 := N-4; \quad N_2 := 2N+1; \quad N_3 := -N-4; \quad N_4 := -2N+1.$$

Let $u_0$ be an $H^s$ function.  We define the iterates $u^{(0)}$, $u^{(1)}$
by
\bas
u^{(0)}_t + \frac{1}{4\pi^2} u^{(0)}_{xxx} &= 0; \quad u^{(0)}(x,0) = u_0(x)\\
u^{(1)}_t + \frac{1}{4\pi^2} u^{(1)}_{xxx} + (u^{(0)})^3 u^{(0)}_x &= 0; \quad u^{(1)}(x,0) = u_0(x).
\end{align*}

In order for the Cauchy problem \eqref{gkdv} to be locally analytically well-posed in $H^s$, it is necessary that the non-linear map $u_0 \mapsto u^{(1)}$ maps $H^s$ to $L^\infty_t H^s_x$, at least for short times $t$ and small $H^s$ norm.  This is because $u^{(1)}$ is the Taylor expansion of $u$ to fourth order in terms of $u_0$ (see e.g. \cite{borg:measures} for further discussion).

We now choose a specific choice of initial data $u_0$, namely
$$ u_0(x) := \eps \sum_{j=0}^4 N_j^{-s} \cos(2\pi N_j x)$$
for some $0 < \eps \ll 1$.  Clearly $u_0$ has $H^s$ norm $O(\eps)$ and mean zero.  The zeroth iterate $u^{(0)}$ is
$$ u^{(0)}(t,x) = \eps \sum_{j=0}^4 N_j^{-s} \cos(2\pi (N_j x - N_j^3 t)).$$
By \eqref{whip}, one can then see that the non-linear expression
$(u^{(0)})^3 u^{(0)}_x$ contains a term of the form
$$ C \eps^4 N_0^{-s} N_1^{-s} N_2^{-s} N_3^{-s} N_4 \cos(2 \pi (N_4 x - (N_4^3 + O(1)) t))$$
for some absolute \emph{non-zero} constant $C$.  From this and a little Fourier analysis one can see that for all non-zero times $|t| \ll 1$, the $N_4$ Fourier coefficient norm of $u^{(1)}(t)$ is $\sim c(t) \eps^4 N^{1-3s}$ for some non-zero quantity $c(t)$ depending only on $t$.  This implies that the $H^s$ norm of $u^{(1)}(t)$ is at least $c(t) N^{1-2s}$.  Since $s < \half$, we obtain $H^s$ analytic ill-posedness by letting $N \to \infty$.  (In fact, by examining the argument more carefully we can show that the solution map is not even $C^4$ in the $H^s$ topology in the $k=3$ case.  For more general $k$, one can show $C^{k+1}$ ill-posedness.).

Note that the above example can be easily modified to also show that Proposition \ref{main-cor} fails for $s < \half$.

\section{Linear Estimates}\label{embedding-sec}

In this section we list some linear embeddings which will be useful in treating non-endpoint cases.

We shall implicitly use the trivial embedding 
\be{triv}
X^{s,b} \subseteq X^{s',b'}
\end{equation}
for $s \geq s'$, $b \geq b'$ frequently in the sequel.
From spatial Sobolev we have 
\be{2-sob}
X^{s,0} = L^2_t H^s_x \subseteq L^2_t L^p_x
\end{equation}
whenever $0 \leq s < \half$ and $2 \leq p \leq \frac{2}{(1-2s)}$, or whenever
$s > \half$ and $2 \leq p \leq \infty$.  Similarly, we have the energy estimate\footnote{We use $a+$ and $a-$ to denote quantities $a+\eps$, $a-\eps$, where $\eps > 0$ is arbitrarily small, and implicit constants are allowed to depend on $\eps$.}
\be{energy}
X^{s,\half+} \subseteq L^\infty_t H^s_x \subseteq L^\infty_t L^p_x
\end{equation}
under the same conditions on $s$ and $p$.  In particular, we have
\be{infty}
X^{\half+,\half+} \subseteq L^\infty_{x,t}.
\end{equation}
By interpolation with the previous estimates we thus have $X^{\half+,\half+} \subseteq L^q_t L^r_x$ for all $2 \leq q,r \leq \infty$.  Interpolating this with \eqref{2-sob} for $s=0$, $p=2$ we obtain
\be{cook}
X^{\half - \delta,\half - \delta} \subseteq L^q_t L^r_x
\end{equation}
for all $0 < \delta < \half$ and $2 \leq q,r < \frac{1}{\delta}$.

From \cite{borg:xsb} we have the Strichartz estimates
\be{strich-4}
X^{0,\frac{1}{3}} \subseteq L^4_{x,t}
\end{equation}
and
$$
X^{0+,\half+} \subseteq L^6_{x,t}.
$$
Interpolating the latter estimate with \eqref{cook} we obtain the improvement
\be{strich-6}
X^{\delta,\half} \subseteq L^q_{x,t}  
\end{equation}
for all $0 < \delta < \half$ and $2 \leq q < \frac{6}{1-2\delta}$. In particular we may take $q = 6$, $6-$, or $6+$.  If we interpolate with \eqref{strich-4} instead we obtain
\be{strich-q}
X^{0+,\half-\sigma} \subseteq L^q_{x,t}
\end{equation}
for all $4 < q < 6$ and $\sigma < 2(\frac{1}{q} - \frac{1}{6})$.

Now we give some embeddings for the $Y^s$ and $Z^s$ spaces.  Since the Fourier transform of an $L^1$ function is continuous and bounded, we have from \eqref{y-def} that
\be{sup-est}
Y^s \subseteq C_t H^s_x \subseteq L^\infty_t H^s_x.
\end{equation}
Let $\eta(t)$ denote a bump function adapted to $[-2,2]$ which equals one on $[-1,1]$.  It is easy to see that multiplication by $\eta(t)$ is a bounded operation on the spaces $Y^s$, $Z^s$, $X^{s,b}$.

Let $S(t)$ denote the evolution operator for the Airy equation:
$$S(t) := \exp(\frac{1}{4\pi^2} t\partial_{xxx}).$$
From the identity
$$ \widehat{\eta(t) S(t) u_0}(\xi,\tau) = \hat u_0(\xi) \hat \eta(\tau - \xi^3)$$
we see that
\be{grow}
\| \eta(t) S(t) u_0 \|_{Y^s} \lesssim \|u_0 \|_{H^s}.
\end{equation}
This homogeneous estimate controls the linear portion of the generalized KdV equation.  To control the Duhamel term we need the following inhomogeneous estimate (cf. \cite{kpv:kdv})
\begin{lemma}\label{invert}
We have
$$ \| \eta(t) \int_0^t S(t-t') F(t')\ dt' \|_{Y^s} \lesssim \| F \|_{Z^s}$$
for any $s$ and arbitrary test functions $F$ on $\T \times \R$.
\end{lemma}

\begin{proof}
Fix $F$; by applying a smooth cutoff one may assume that $F$ is supported on $\T \times [-3, 3]$. 

Let $a(t) := \sgn(t) \tilde \eta(t)$, where $\tilde \eta$ is a bump function on $[-10, 10]$ which equals 1 on $[-5, 5]$.  From the identity
$$ \chi_{[0,t]}(t') = \frac{1}{2} (a(t') - a(t-t'))$$
for all $t \in [-2,2]$ and $t' \in [-3,3]$, we see that we may write $\eta(t) \int_0^t S(t-t') F(t')\ dt'$ as a linear combination of
\be{t1}
\eta(t) S(t) \int_{\R} a(t') S(-t') F(t')\ dt'
\end{equation}
and
\be{t2}
\eta(t) \int_{\R} a(t-t') S(t-t') F(t')\ dt'.
\end{equation}
Consider the contribution of \eqref{t1}.  By \eqref{grow} it suffices to show that
$$ 
\left\| \int a(t') S(-t') F(t')\ dt' \right\|_{H^s} \lesssim \| F\|_{Z^s}.$$
Observe that the Fourier transform of $\int a(t') S(-t') F(t')\ dt'$ at $\xi$ is given by
$$ \int \hat a(\tau - \xi^3) \hat F(\xi, \tau)\ d\tau.$$
Since one has the easily verified bound
\be{a-bound}
\hat a(\lambda) = O( \langle \lambda \rangle^{-1}),
\end{equation} 
the claim then follows from \eqref{z-def}.

Now consider the contribution of \eqref{t2}.
We may discard the $\eta(t)$ cutoff.  The spacetime Fourier transform of $\int_{\R} a(t-t') S(t-t') F(t')\ dt'$ at $(\tau, \xi)$ is equal to $\hat a(\tau - \xi^3) \hat F(\tau,\xi)$.  The claim then follows from \eqref{a-bound}, \eqref{y-def}, \eqref{z-def}.  
\end{proof}

Finally, we shall need the following duality relationship between $Y^s$ and $Z^{-s}$.

\begin{lemma}\label{duality}
We have
$$ \left| \int\int \chi_{[0,1]}(t) u(x,t) v(x,t)\ dx dt \right| \lesssim \| u \|_{Y^s} \| v \|_{Z^{-s}}$$
for any $s$ and any $u$, $v$ on $\T \times \R$.
\end{lemma}

\begin{proof}
Without the cutoff $\chi_{[0,1]}$ this would be an immediate consequence of the duality of $X^{s,\half}$ and $X^{-s,-\half}$.  However, these spaces are not preserved by rough cutoffs, and one requires a little more care.

By writing $\chi_{[0,1]}$ as the difference of two signum functions, it suffices to show that
$$ |\int\int \sgn(t) u v\ dx dt| \lesssim \| u \|_{Y^s} \| v \|_{Z^{-s}}.$$
By Plancherel we can write the left-hand side as
$$ C |p.v. \int\int\int \hat u(\xi, \tau) \hat v(\xi, \tau') \frac{d \xi d\tau d\tau'}{\tau - \tau'}|.$$

Partition $u = \sum_{j = 0}^\infty u_j$, where $u_j$ has Fourier support on the region $\langle \xi - \tau^3 \rangle \sim 2^j$.  Similarly partition $v = \sum_{k=0}^\infty v_k$.  We can thus estimate the above by
\be{combined}
C \sum_{j, k \geq 0} |p.v. \int\int\int \langle \xi \rangle^s \hat u_j(\xi, \tau) \langle \xi \rangle^{-s} \hat v_k(\xi, \tau') \frac{d \xi d\tau d\tau'}{\tau - \tau'}|.
\end{equation}
First consider the contribution of the case $|j-k| \lesssim 1$.  In this case we use the $L^2$ boundedness of the Hilbert transform and Cauchy-Schwarz to estimate the above by
$$ C \sum_{j, k \geq 0: |j-k| \lesssim 1} \| \langle \xi \rangle^s \hat u_j \|_2
\| \langle \xi \rangle^{-s} \hat v_k \|_2.$$
Since $|j-k| \lesssim 1$, we may estimate this by
$$ C \sum_{j, k \geq 0: |j-k| \lesssim 1} \| \langle \xi \rangle^s \langle \xi - \tau^3 \rangle^{\half} \hat u_j \|_2
\| \langle \xi \rangle^{-s} \langle \xi - \tau^3 \rangle^{-\half} \hat v_k \|_2,$$
which by another Cauchy-Schwarz is bounded by
$$ C \| u \|_{s,\half} \|v\|_{-s,-\half} \lesssim \| u \|_{Y^s} \| v \|_{Z^s}$$
as desired.

Now consider the contribution when $|j-k| \gg 1$.  In this case we observe that $|\tau - \tau'| \gtrsim \langle \tau' - \xi^3 \rangle$, so we may estimate \eqref{combined} by
$$
C
p.v. \int\int\int \langle \xi \rangle^s |\hat u(\xi, \tau)|
\frac{\langle \xi \rangle^{-s} |\hat v(\xi, \tau')|}{\langle \tau' - \xi^3 \rangle} d \xi d\tau d\tau',$$
which by Fubini and Cauchy-Schwarz is bounded by
$$ C \| \langle \xi \rangle^s \hat u \|_{L^2_\xi L^1_\tau} \| 
\frac{\langle \xi \rangle^{-s} \hat v}{\langle \tau' - \xi^3 \rangle} \|_{L^2_\xi L^1_{\tau'}} \lesssim \| u \|_{Y^s} \|v\|_{Z^s}$$
as desired.
\end{proof}

\section{Reduction to a multiplier bound}\label{Non-endpoint-sec}

In this section we make some preliminary reductions for Theorem \ref{main}, exploiting the ``denominator games'' of Bourgain and reducing matters to a multilinear multiplier estimate.

Fix $s$, $u_i$.  Since the spaces $Y^s$ and $X^{s-1,\half}$ are defined using the size of the Fourier transform, we may assume that the Fourier transforms of the $u_i$ are non-negative.

The desired estimate \eqref{main-est} is trivial for $k=1$, so we may assume $k \geq 2$.

The Fourier transform of $\prod_{i=1}^k u_i$ is given by
$$ \widehat{\prod_{i=1}^k u_i}(\xi,\tau) = \int_* \prod_{i=1}^k \hat u_i(\xi_i, \tau_i)$$
where $\int_*$ denotes an integration over the set where $\xi = \xi_1 + \ldots + \xi_k$, $\tau = \tau_1 + \ldots + \tau_k$.

First consider the contribution where 
\be{tau1-dom}
\langle \tau - \xi^3 \rangle \lesssim \langle \tau_1 - \xi_1^3 \rangle.
\end{equation}
In this case it suffices to show that
$$ \| \prod_{i=1}^k u_i \|_{s-1,0} \lesssim \| u_1 \|_{s,0}
\prod_{i=2}^k \| u_i \|_{Y^{s}}.$$
On the other hand, from Sobolev and H\"older, then another Sobolev and \eqref{sup-est} we have
\bas
\|  \prod_{i=1}^k u_i \|_{s-1,0} &\lesssim
\|  \prod_{i=1}^k u_i \|_{L^2_t L^{1+}_x}\\
&\lesssim \| u_1 \|_{L^2_t L^{k+}_x} \prod_{i=2}^k \| u_i \|_{L^\infty_t L^{k+}_x}\\
&\lesssim \| u_1 \|_{L^2_t H^s_x} \prod_{i=2}^k \| u_i \|_{L^\infty_t H^s_x}\\
&\lesssim \| u_1 \|_{s,0} \prod_{i=2}^k \| u_i \|_{Y^s}
\end{align*}
as desired.

From the above and symmetry, we may thus assume that
\be{tau-dom}
\langle \tau - \xi^3  \rangle \gg \langle \tau_i - \xi_i^3 \rangle
\end{equation}
for all $1 \leq i \leq k$.  In particular we have
$$ 1 \ll \langle \tau - \xi^3  \rangle \sim |\xi^3 - \sum_{i=1}^k \xi_i^3|.$$
It thus suffices to show that
\be{key-targ}
\|
\int_* \langle \xi\rangle^{s-1} |\xi^3 - \sum_{i=1}^k \xi_i^3|^{\half} \prod_{i=1}^k \hat u_i(\xi_i, \tau_i) \|_{L^2_\tau L^2_\xi}
\lesssim \prod_{i=1}^k \| u_i \|_{s,\half}.
\end{equation}
We now rewrite \eqref{key-targ} using the notation of \cite{tao:xsb}.  For all $n \geq 2$ and all symbols $m(\xi_1, \tau_1, \ldots, \xi_n, \tau_n)$ defined on the region
$$ \Gamma_n := \{ (\xi_1, \ldots, \tau_n) \in (\Z \times \R)^n: \xi_1 + \ldots + \xi_n = \tau_1 + \ldots + \tau_n = 0 \},$$
define the norm
$$ \| m \|_{[n; \Z \times \R]}$$
to be the best constant such that one has the bound
$$ \left| \int_{\Gamma_n} m(\xi_1, \ldots, \tau_n) \prod_{i=1}^n f_i(\xi_i, \tau_i) \right|
\leq \| m \|_{[n; \Z \times \R]} \prod_{i=1}^n \| f_i \|_{L^2_{\xi_i} L^2_{\tau_i}}.$$
By duality \eqref{key-targ} can now be written as
$$
\| \frac{|\sum_{i=1}^{k+1} \xi_i^3|^{\half}}{\langle \xi_{k+1} \rangle^{1-s}
\prod_{i=1}^k \langle \xi_i \rangle^{s} \langle \tau_i - \xi_i^3  \rangle^{\half}}
\|_{[k+1; \Z \times \R]} \lesssim 1.$$
We would like to use the numerator to cancel some of the denominator.  Our tool for doing this is 
\begin{lemma}\label{cubic-sum}
If $|\xi_1| \geq \ldots \geq |\xi_{k+1}|$ and $\sum_{i=1}^{k+1} \xi_i = 0$, then
$$ \sum_{i=1}^{k+1} \xi_i^3 = O(|\xi_1| |\xi_2| |\xi_3|).$$
\end{lemma}

\begin{proof}
From the estimate
$$ \xi_1^3 + \xi_2^3 = (\xi_1 + \xi_2)(\xi_1^2 - \xi_1 \xi_2 + \xi_2^2) = O(|\xi_1 + \xi_2| |\xi_1|^2)$$
we see that
$$ \sum_{i=1}^{k+1} \xi_i^3 = O(|\xi_1 + \xi_2| |\xi_1|^2) + O(|\xi_3|^3).$$
Since $|\xi_1 + \xi_2| = O(|\xi_3|)$ and $|\xi_3| \leq |\xi_2| \sim |\xi_1|$, the claim follows.
\end{proof}
From Lemma \ref{cubic-sum} and symmetry of the $\xi_1, \ldots, \xi_k$ variables, it thus suffices to show the estimates
\be{est-2}
\| \frac{|\xi_1|^{\half} |\xi_2|^{\half} |\xi_{k+1}|^{\half}}{
\langle \xi_{k+1} \rangle^{1-s}
\prod_{i=1}^k \langle \xi_i \rangle^s \lambda_i^{\half}}
\|_{[k+1; \Z \times \R]} \lesssim 1
\end{equation}
and (when $k \geq 3$)
\be{est-1}
\| \frac{|\xi_1|^{\half} |\xi_2|^{\half} |\xi_3|^{\half}}{
\langle \xi_{k+1} \rangle^{1-s}
\prod_{i=1}^k \langle \xi_i \rangle^s \lambda_i^{\half}}
\|_{[k+1; \Z \times \R]} \lesssim 1
\end{equation}
where we adopt the notation
$$ \lambda_i := \langle \tau_i - \xi_i^3 \rangle.$$

In the next section we prove these estimates in the non-endpoint case $s > \half$.  Then in the following two sections we resolve the more difficult endpoint case $s=\frac{1}{2}$.

\section{The non-endpoint case}\label{1-sec}

We now prove \eqref{est-2}, \eqref{est-1} in the non-endpoint case $s > \half$.

First consider \eqref{est-2}.  By symmetry and the Comparison Principle (\cite{tao:xsb}, Lemma 3.1) we may assume that $|\xi_1| \geq \ldots \geq |\xi_k|$.  In particular we have $|\xi_{k+1}| \lesssim |\xi_1|$, so that we have the pointwise inequality
$$
\frac{ |\xi_1|^{\half} |\xi_2|^{\half} |\xi_{k+1}|^{\half} }{\langle \xi_1 \rangle^s
\langle \xi_2 \rangle^s \langle \xi_{k+1} \rangle^{1-s}}
\lesssim \frac{1}{\langle \xi_2 \rangle^{s-\half}}.$$
To show \eqref{est-2} it thus suffices to show that
$$
\| \frac{1}{ \langle \xi_2 \rangle^{s-\half}
(\prod_{i=3}^{k} \langle \xi_i \rangle^s)
\prod_{i=1}^k \lambda_i^{\half}}
\|_{[k+1; \Z \times \R]} \lesssim 1
$$ 
which by duality becomes
$$
\| u_1 \ldots u_k \|_{L^2_{x,t}} \lesssim
\| u_1 \|_{0,\half} \| 
u_2 \|_{s-\half,\half} \prod_{i=3}^k \| u_i \|_{s,\half}.
$$
However, by H\"older we may take $u_1$ in $L^4_{x,t}$, $u_2$ in $L^6_{x,t}$, and the other $u_i$ in $L^{12(k-2)}_{x,t}$.  The claim then follows from \eqref{strich-4}, \eqref{strich-6}, and \eqref{cook}.

Now consider \eqref{est-1}.  By duality as before it thus suffices to show that
$$\| u_1 \ldots u_k \|_{L^2_t H^{s-1}_x} \lesssim
\| u_1 \|_{s-\half,\half} \| u_2 \|_{s-\half,\half} \| u_3 \|_{s-\half,\half}
\prod_{i=4}^k \| u_i \|_{s,\half}.$$
We may of course estimate the $L^2_t H^{s-1}_x$ norm by the $L^2_{x,t}$ norm.  When $k \geq 4$ we use H\"older to take $u_1$, $u_2$, $u_3$ in $L^{6+}_{x,t}$ and the other $u_i$ in $L^{\infty-}_{x,t}$, then use \eqref{strich-6} and \eqref{cook}.  When $k=3$ we just take $u_1$, $u_2$, $u_3$ in $L^6_{x,t}$ and use \eqref{strich-6}.  This completes the proof of Theorem \ref{main} in the non-endpoint case.

\section{Some elementary number theory}\label{overview-sec}

In the non-endpoint arguments in the preceding section we used the $L^6$ Strichartz estimate \eqref{strich-6} from \cite{borg:xsb}.  Let us quickly review the method of proof for this estimate.  It suffices to show that
\be{trilinear}
\| u_1 u_2 u_3 \|_{L^2_{x,t}} \lesssim \| u_1 \|_{0+,\half+} \| u_2 \|_{0+,\half+} \| u_3 \|_{0+,\half+}.
\end{equation}
By standard Cauchy-Schwarz arguments (see \cite{borg:xsb}, or apply \cite{tao:xsb} Proposition 5.1 and Lemma 3.9) this result would obtain if we knew that the number of integer solutions of size $O(N)$ to the equations
$$ \xi = \xi_1 + \xi_2 + \xi_3; \tau = \xi_1^3 + \xi_2^3 + \xi_3^3$$
was bounded by $O(N^{0+})$ uniformly for non-zero $(\xi,\tau)$ and $N \geq 1$.  But this follows from the identity (cf. \cite{borg:xsb})
$$ \tau - \xi^3 = 3(\xi_1 + \xi_2)(\xi_2 + \xi_3)(\xi_3 + \xi_1),$$
and the well-known observation (see e.g. \cite{apostol})
\be{divisors}
\hbox{Every non-zero integer } \lambda \hbox{ has at most } O(|\lambda|^{0+}) \hbox{ factors}
\end{equation}
from elementary number theory.  (The case $\tau - \xi^3 = 0$ needs to be dealt with separately).

It would be very convenient if one could replace $0+$ by $0$ in the above arguments; indeed, the endpoint estimate would then follow by a variant of the preceding arguments.   Even with the epsilon loss in exponents, \eqref{strich-6} is still strong enough to treat a large portion of \eqref{est-2} and \eqref{est-1} in the endpoint case.

We do not know how to prove the endpoint of \eqref{strich-6} directly (see \cite{borg:xsb} for some further discussion of this issue and of a related $L^8_{x,t}$ conjecture).  However, we will be able to remove the epsilon in \eqref{trilinear} when $u_3$ (for instance) has much smaller frequency than the other two functions, and this is enough to treat the remaining cases for \eqref{est-2} and \eqref{est-1}.  

To achieve this we shall rely on the following variant of \eqref{divisors}.

Let $\xi$, $\lambda$, $N$, $L$ be integers such that
$0 < L, N \ll |\xi| \lesssim |\lambda|$.  We consider the quantity
\be{interval}
\# \{ (l,n) \in \Z^2: |l-\lambda| \lesssim L; |n-\xi| \lesssim N; n \left| l \right. \},
\end{equation}
where we use the notation $a | b$ to denote that $a$ divides $b$.  In other words, for all $l = \lambda + O(L)$, we count the divisors of $l$ which lie in the interval $\xi + O(N)$.  From \eqref{divisors} we may clearly bound \eqref{interval} by $|\lambda|^{0+} L$.  The purpose of the following lemma is to remove the $|\lambda|^{0+}$ under some additional assumptions.

\begin{lemma}[Few divisors in small intervals]\label{div} 
Let $\xi$, $\lambda$, $N$, $L$ be as above.  Then 
$$ \eqref{interval} \leq N.$$
If we make the further assumptions
$$ |\lambda| \lesssim |\xi|^3; \quad 0 < N \ll |\xi|^{\frac{1}{6}}$$
then we have the additional bound
$$ \eqref{interval} \leq 3L.$$
\end{lemma}

The assumptions and bounds are far from best possible, but suffice for our purposes.

\begin{proof}
Since $L \ll |\xi|$, we have
$$ \# \{ l \in \Z: |l-\lambda| \lesssim L; n \left| l \right. \} \leq 1$$
for all $|n| \sim |\xi|$, and the first bound follows.  

Now we prove the second bound.  It suffices to show that
$$ \# \{ n \in \Z: |n-\xi| \lesssim N; n \left| l \right. \} \leq 3$$
for all $|l| \sim |\xi|^3$.

Fix $l$, and suppose for contradiction that there were four integers $n_1, n_2, n_3, n_4$ in the above set.  Then from elementary number theory we see that
$$ \prod_{1 \leq i \leq 4} n_i \left| l \prod_{1 \leq i < j \leq 4} \gcd(n_i,n_j)\right. .$$
On the other hand, each $n_i$ has size $\sim |\xi|$, and from the Euclidean algorithm we see that $\gcd(n_i,n_j) \lesssim N \ll |\xi|^{\frac{1}{6}}$.  Since $l = O(|\xi|^3)$, we obtain the desired contradiction. 
\end{proof}

From this lemma and the previous Cauchy-Schwarz argument we can obtain various partial endpoint versions of \eqref{trilinear}.  For instance, we can prove
$$
\| u_1 u_2 u_3 \|_{L^2_{x,t}} \lesssim \| u_1 \|_{0,\half+} \| u_2 \|_{0,\half+} \| u_3 \|_{0,\half+}$$
when $u_j$ has Fourier support in the region $N_j$ and $N_1 \gg N_2 \gg N_3$, under the additional assumption\footnote{One can view this as a kind of ``trilinear improved Strichartz inequality'', similar to the bilinear improvements to Strichartz inequalities in e.g. \cite{borg:book}.} $N_3 \ll N_1^{\frac{1}{3}}$.  We sketch the argument as follows.  We repeat the proof of \eqref{trilinear} but observe that we may assume $\xi = O(N_1)$, $\xi_1 + \xi_2 = -\xi - \xi_3 = -\xi + O(N_3)$, and $\tau - \xi_3^3 = O(N_1^3)$.  If we discard the special case $\tau - \xi_3^3 = 0$, we may invoke the second part of Lemma \ref{div} (with $L=1$) and conclude that for fixed $\tau$, $\xi$, there are at most 3 values of $\xi_1 + \xi_2$.  Fixing $\xi_1 + \xi_2$ determines $\xi_3$, and hence $\xi_1^3 + \xi_2^3$, which determines $\xi_1$ and $\xi_2$ up to permutations.   So we have at most 6 integer solutions $(\xi_1,\xi_2,\xi_3)$ (rather than $O(N_1^{0+})$), and the claim follows.  

We will not use this partial endpoint result directly in the sequel, but arguments with the above flavor will be used to compensate for the $0+$ loss in the $L^6$ Strichartz estimate at various junctures.

\section{The proof of \eqref{est-2} in the endpoint case.}\label{2-prf}

We now prove \eqref{est-2} in the endpoint case $s = \half$. 

By duality, it suffices to show that
$$ \| u_1 \ldots u_k \|_{L^2_{x,t}} \lesssim 
\| u_1\|_{0,\half} \| u_2 \|_{0,\half} \prod_{i=3}^k \| u_i \|_{\half,\half}.$$

We will actually prove the stronger estimate 
\be{est2-plus}
\| u_1 \ldots u_k \|_{L^2_{x,t}} \lesssim 
\| u_1\|_{0,\half-\delta} \| u_2 \|_{0,\half-\delta} \prod_{i=3}^k \| u_i 
\|_{\half - \delta,\half -\delta}
\end{equation}

for some small $0 < \delta \ll 1$; this improved estimate shall be useful for the large data applications.  We may assume that $k \geq 3$ since the claim follows from \eqref{strich-4} and H\"older otherwise.  

By symmetry and interpolation\footnote{Specifically, we permute $u_3, \ldots, u_k$ and perform multilinear complex interpolation to obtain the ``centroid'' of  all the permuted estimates.  The small losses of $0+$ in some of the indices  are more than compensated for by the gains of $\frac{1}{100}$ in other indices, and so in the final estimate we will have some sort of gain $\delta > 0$ throughout.}, it suffices to show 

$$ \| u_1 \ldots u_k \|_{L^2_{x,t}} \lesssim 
\| u_1\|_{0, \frac{1}{2}-\frac{1}{100}} \| u_2 \|_{0, \frac{1}{2}-\frac{1}{100}} \| u_3 \|_{\frac{1}{2}-\frac{1}{100},\frac{1}{2}-\frac{1}{100}}
\prod_{i=4}^k \| u_i \|_{\half+,\half+}.$$

The exponent $\frac{1}{100}$ has no special significance, and could be replaced by any other small constant.

By \eqref{infty}, we need only show the trilinear estimate
$$ \| u_1 u_2 u_3\|_{L^2_{x,t}} \lesssim 
\| u_1\|_{0, \frac{1}{2}-\frac{1}{100}} \| u_2 \|_{0, \frac{1}{2}-\frac{1}{100}} \| u_3 \|_{\frac{1}{2}-\frac{1}{100},\frac{1}{2}-\frac{1}{100}}.$$
By dyadic decomposition, it suffices to show that
$$ \| u_1 u_2 u_3\|_{L^2_{x,t}} \lesssim N_3^{\frac{1}{2}-\frac{1}{100}-}
\| u_1\|_{0, \frac{1}{2}-\frac{1}{100}} \| u_2 \|_{0, \frac{1}{2}-\frac{1}{100}} \| u_3 \|_{0,\frac{1}{2}-\frac{1}{100}}$$
for all $N_3 \geq 1$, where $u_3$ is supported on $\langle \xi_3 \rangle \sim N_3$.

Fix $N_3$. By duality we reduce to
$$ \| \frac{\chi_{|\xi_3| \sim N_3}
}{
(\lambda_1 \lambda_2 \lambda_3)^{\frac{1}{2}-\frac{1}{100}}
} \|_{[4; \Z \times \R]} \lesssim N_3^{\frac{1}{2}-\frac{1}{100}-},$$
where we have again adopted the notation
$$ \lambda_i := \langle \tau_i - \xi_i^3 \rangle.$$
By another dyadic decomposition it suffices to show that
\be{eft-targ}
\| \chi_{|\xi_3| \sim N_3 }
\chi_{\lambda_1 \sim L_1}
\chi_{\lambda_2 \sim L_2}
\chi_{\lambda_3 \sim L_3}
\|_{[4; \Z \times \R]} \lesssim 
N_3^{\frac{1}{2}-\frac{1}{100}-} (L_1 L_2 L_3)^{\frac{1}{2}-\frac{1}{100}-}
\end{equation}
for all $L_1,L_2,L_3 \gtrsim 1$.

Fix $L_1,L_2,L_3$.  By the Conjugation Lemma (\cite{tao:xsb}, Corollary 3.8, or the observation that $\| uv \|_2 = \|u \overline{v} \|_2$) we may assume that $\xi_1$, $\xi_2$ have the same sign.  By symmetry we may thus take $\xi_1$, $\xi_2$ non-negative.

Let us first deal with the case\footnote{This is of course a rather large portion of the integral.  The reason we can handle so much of the integral so easily is because we are at the endpoint of a much easier non-endpoint result, and the non-endpoint arguments contain some ``slack'' in the ``$b$'' index of the $X^{s,b}$ norm.  In this case we can borrow some of this slack to create a small amount of room in the ``$s$'' index, at which point we are no longer constrained by our failure to prove the endpoint $L^6_{x,t}$ Strichartz estimate.} when $\xi_2 \lesssim (N_3 L_1 L_2 L_3)^{10}$.  In this case we may borrow some regularity from $u_3$ to place on $u_2$, and it suffices by duality to show that 

$$ \|u_1 u_2 u_3\|_2 \lesssim \| u_1 \|_{0, \frac{1}{2}-\frac{2}{100}} \| u_2 \|_{0+,\frac{1}{2}-\frac{2}{100}}
\| u_3 \|_{\frac{1}{2}-\frac{2}{100},\frac{1}{2}-\frac{2}{100}}.$$
But this follows by (for instance) taking $u_1$ in $L^4_{x,t}$, $u_2$ in $L^5_{x,t}$, and $u_3$ in $L^{20}_{x,t}$, and using \eqref{strich-4}, \eqref{strich-q}, \eqref{cook}.
From this and symmetry we may restrict ourselves to the case $\xi_1, \xi_2 \gg (N_3 L_1 L_2 L_3)^{10}$.

By Cauchy-Schwarz (\cite{tao:xsb}, Lemma 3.9) it suffices to show that
$$ \int \chi_{\xi_1, \xi_2 \gg (N_3 L_1 L_2 L_3)^{10}} \chi_{|\xi_3| \sim N_3 }
\chi_{\lambda_1 \sim L_1}
\chi_{\lambda_2 \sim L_2}
\chi_{\lambda_3 \sim L_3}\ d\xi_1 d\xi_2 d\tau_1 d\tau_2
\lesssim 
N_3^{1-\frac{2}{100}-} (L_1 L_2 L_3)^{1-\frac{2}{100}-}$$
for all $\xi_4$, $\tau_4$, where $\xi_3$, $\tau_3$ is given by $\xi_1 + \xi_2 + \xi_3 + \xi_4 = \tau_1 + \tau_2 + \tau_3 + \tau_4 = 0$.

Fix $\xi_4$, $\tau_4$.  Performing the $\tau$ integrals\footnote{Because of the restriction $\tau_1 + \tau_2 + \tau_3 + \tau_4 = 0$ we have some choice in which two $\tau_j$ indices to integrate.  To obtain the optimal factor $L_{min} L_{med}$ one should integrate the two $\tau_j$ corresponding to the quantities $L_{min}$, $L_{med}$.}, we reduce to
$$ L_{min} L_{med} \# \{ (\xi_1, \xi_2, \xi_3) \in \Z^3: \xi_1 + \xi_2 + \xi_3 + \xi_4 = 0;
\xi_1, \xi_2 \gg (N_3 L_1 L_2 L_3)^{10}; |\xi_3| \sim N_3;$$
$$ |\xi_1^3 + \xi_2^3 + \xi_3^3 + \tau_4| \lesssim L_{max} \}
\lesssim
N_3^{1-\frac{2}{100}-} (L_1 L_2 L_3)^{1-\frac{2}{100}-}$$
where $L_{min} \leq L_{med} \leq L_{max}$ are the minimum, median, and maximum of $L_1, L_2, L_3$ respectively.  It will suffice to show that
$$ \# \{ (\xi_1, \xi_2, \xi_3) \in \Z^3: \xi_1 + \xi_2 + \xi_3 + \xi_4 = 0;
\xi_1, \xi_2 \gg (N_3 L_{max})^{10}; |\xi_3| \sim N_3;$$
$$ |\xi_1^3 + \xi_2^3 + \xi_3^3 + \tau_4| \lesssim L_{max} \}
\lesssim
N_3^{1-\frac{2}{100}-} (L_{max})^{1-\frac{6}{100}-}$$
Since $\xi_1 + \xi_2 + \xi_3 + \xi_4 = 0$, we have the identity
\be{4-way}
\xi_1^3 + \xi_2^3 + \xi_3^3 + \xi_4^3 = 3(\xi_1 + \xi_2)(\xi_2 + \xi_3)(\xi_3 + \xi_1)
\end{equation}
and so we reduce to showing
\be{card}
\begin{split}
\# \{ (\xi_1, \xi_2, \xi_3) \in \Z^3: &\xi_1 + \xi_2 + \xi_3 + \xi_4 = 0;
\xi_1, \xi_2 \gg (N_3 L_{max})^{10}; |\xi_3| \sim N_3; \\
&(\xi_1 + \xi_2) (\xi_2 + \xi_3) (\xi_3 + \xi_1) \in I \}
\lesssim N_3^{1-\frac{2}{100}-} L_{max}^{1-\frac{6}{100}-}
\end{split}
\end{equation}
where $I$ is the set
$$ I := \{ l \in \Z: |l + (\tau_4 - \xi_4^3 )/3| \lesssim L_{max} \}.$$
Since $\xi_1$, $\xi_2$ are positive and much larger than $\xi_3$, we see that $|\xi_4| \gtrsim   |\xi_1|, |\xi_2|$, and that
$$ |(\xi_1 + \xi_2) (\xi_2 + \xi_3) (\xi_3 + \xi_1)| \lesssim |\xi_1 + \xi_2|^3
\lesssim |\xi_4|^3.$$
We may therefore assume that $|\tau_4 - \xi_4^3| \lesssim |\xi_4|^3$, since \eqref{card} vanishes otherwise.

From Lemma \ref{div} thus we see that
$$ \# \{ (n,l): |n + \xi_4| \lesssim N_3; l \in I; n \left| l \right. \} \lesssim
\min(N_3, L_{max}).$$
From this and the previous discussion we see that there are at most $O(\min(N_3,L_{max}))$ possible values of $\xi_1 + \xi_2$ which can contribute to \eqref{card}.  But from elementary algebra we see that each value of $\xi_1 + \xi_2$ contributes at most $O(1)$ elements to \eqref{card}.  The claim then follows.

\section{The proof of \eqref{est-1} in the endpoint case.}\label{1-prf}

It remains to prove \eqref{est-1}, which we rewrite as
$$ |\int u_1 \ldots u_{k+1}\ dx dt| \lesssim
\| u_1 \|_{0,\half}
\| u_2 \|_{0,\half}
\| u_3 \|_{0,\half}
\| u_{k+1} \|_{\half,0}
\prod_{i=4}^k \| u_i \|_{\half,\half}.$$
We shall actually prove the stronger estimate
\be{est1-plus}
|\int u_1 \ldots u_{k+1}\ dx dt| \lesssim
\| u_1 \|_{0,\half-\delta}
\| u_2 \|_{0,\half-\delta}
\| u_3 \|_{0,\half-\delta}
\| u_{k+1} \|_{\half-\delta,0}
\prod_{i=4}^k \| u_i \|_{\half-\delta,\half - \delta}.
\end{equation}
for some $0 < \delta \ll 1$.  Again, this improved estimate shall be useful for large data applications.

We shall prove \eqref{est1-plus} in the case $k \geq 4$.  The $k=3$ case is slightly simpler, and can be obtained by a routine modification of the following argument.

We first observe that
$$ |\int u_1 \ldots u_{k+1}\ dx dt| \lesssim
\| u_1 \|_{0,\half - \frac{1}{6}}
\| u_2 \|_{0,\half - \frac{1}{6}}
\| u_3 \|_{0,\half+}
\| u_{k+1} \|_{\half+,0}
\prod_{i=4}^k \| u_i \|_{\half+,\half+}.$$
Indeed, this follows by taking $u_1$, $u_2$ in $L^4_{x,t}$, $u_3$ in $L^\infty_t L^2_x$,
$u_{k+1}$ in $L^2_t L^\infty_x$, and all other $u_i$ in $L^\infty_{x,t}$, and then using \eqref{strich-4}, \eqref{energy}, \eqref{2-sob}, and \eqref{infty}.

By symmetry and interpolation\footnote{Specifically, we permute $\{u_1, u_2, u_3\}$ and $\{u_4, \ldots, u_{k+1}\}$ independently in this estimate and the previous one, and then use multilinear complex interpolation (see e.g. \cite{bergh:interp}) to obtain the centroid of all the permuted estimates.  The small losses of $0+$ in some indices will be more than compensated for by the gains of $\frac{1}{100}$ and $\frac{1}{6}$ in other indices, so in the interpolated estimate one will have some non-zero gain $\delta > 0$ throughout.} it thus suffices to show that
$$ |\int u_1 \ldots u_{k+1}\ dx dt| \lesssim
\| u_1 \|_{0,\half+}
\| u_2 \|_{0,\half+}
\| u_3 \|_{0,\half+}
\| u_{k+1} \|_{\frac{1}{2}-\frac{1}{100},0}
\| u_4 \|_{\frac{1}{2}-\frac{1}{100},\frac{1}{2}-\frac{1}{100}}
\prod_{i=5}^k \| u_i \|_{\half+,\half+}.$$

From the fractional Leibniz rule and Sobolev (or see e.g. \cite{tao:xsb}, Corollary 3.16) we have
$$ \left\| (\prod_{i=5}^k f_i) f_{k+1} \right\|_{H^{\frac{1}{2} - \frac{1}{100}}}
\lesssim \left\| f_{k+1} \right\|_{H^{\frac{1}{2} - \frac{1}{100}}} \prod_{i=5}^k
\| f_i \|_{H^{\frac{1}{2}+}}$$
and so by \eqref{energy} we have
$$ \| (\prod_{i=5}^k u_i) u_{k+1} \|_{\frac{1}{2} - \frac{1}{100},0}
\lesssim \| u_{k+1} \|_{\frac{1}{2} - \frac{1}{100},0} 
\prod_{i=5}^k \| u_i \|_{\half+,\half+}.$$
We thus reduce to the quintilinear estimate
$$ |\int u_1 u_2 u_3 u_4 u_5\ dx dt| \lesssim
\| u_1 \|_{0,\half+}
\| u_2 \|_{0,\half+}
\| u_3 \|_{0,\half+}
\| u_4 \|_{\frac{1}{2}-\frac{1}{100},\frac{1}{2}-\frac{1}{100}}
\| u_5 \|_{\frac{1}{2}-\frac{1}{100},0}.$$
By dyadic decomposition it suffices to show that
$$ |\int u_1 u_2 u_3 u_4 u_5\ dx dt| \lesssim
(N_4 N_5 L_4)^{\frac{1}{2}-\frac{1}{100}}
\| u_1 \|_{0,\half+}
\| u_2 \|_{0,\half+}
\| u_3 \|_{0,\half+}
\| u_4 \|_{0,0}
\| u_5 \|_{0,0}$$
for all $N_4, L_4, N_5 \geq 1$, where $u_4$ and $u_5$ have Fourier support in the regions $\langle \xi_4 \rangle \sim N_4$, $\lambda_4 \sim L_4$ and $\langle \xi_5 \rangle \sim N_5$ respectively.

Let us first consider the contribution where $u_3$ is supported in the region
$\langle \xi_3 \rangle \lesssim (N_4 N_5 L_4)^{10}$.  In this case it suffices to show that
$$ |\int u_1 u_2 u_3 u_4 u_5\ dx dt| \lesssim
\| u_1 \|_{0,\half+}
\| u_2 \|_{0,\half+}
\| u_3 \|_{0+,\half+}
\| u_4 \|_{\frac{1}{2}-\frac{2}{100},\frac{1}{2}-\frac{2}{100}}
\| u_5 \|_{\frac{1}{2}-\frac{2}{100},0},$$
since we can borrow some powers of $N_4 N_5 L_4$ to yield a little regularity on $u_3$.  But this estimate can be achieved by taking (for instance) $u_1$ in $L^4_{t,x}$, $u_2$ in $L^\infty_t L^2_x$, $u_3$ in $L^6_{t,x}$, $u_4$ in $L^{12}_t L^{24}_x$, and $u_5$ in $L^2_t L^{24}_x$, and applying \eqref{strich-4}, \eqref{energy}, \eqref{strich-6}, \eqref{cook}, \eqref{2-sob}.

Thus we may assume that $u_3$ is supported in the region $\langle \xi_3 \rangle \gg (N_4 N_5 L_4)^{10}$, and similarly for $u_1$, $u_2$.

By averaging arguments (\cite{tao:xsb}, Proposition 5.1) we may thus assume that $u_1$, $u_2$, $u_3$ have Fourier support on the region
$$\Omega := \{ (\xi, \tau) \in \Z \times \R: |\xi| \gg (N_4 N_5 L_4)^{10}; \tau = \xi^3 + O(1) \}.$$
We may thus rewrite our desired estimate as
$$
\| 
\chi_{\Omega}(\xi_1, \tau_1)
\chi_{\Omega}(\xi_2, \tau_2)
\chi_{\Omega}(\xi_3, \tau_3)
\chi_{|\xi_4| \lesssim N_4}
\chi_{|\xi_5| \lesssim N_5}
\chi_{\lambda_4 \lesssim L_4}
\|_{[5; \Z \times \R]}
\lesssim
(N_4 N_5 L_4)^{\frac{1}{2}-\frac{1}{100}}.$$

This quintilinear estimate is too complex to handle directly.  The strategy will be to reduce this estimate to a quartilinear estimate, and (in some cases) further to a trilinear estimate.

By symmetry we may assume that $|\xi_1| \leq |\xi_2| \leq |\xi_3|$; since $|\xi_4| + |\xi_5| \lesssim (N_4 N_5 L_4)^{10} \ll |\xi_1|, |\xi_2|, |\xi_3|$, we thus see that $|\xi_2| \sim |\xi_1 + \xi_2|$.  Thus we may freely insert a factor of
$$ \chi_\Sigma(\xi_1 + \xi_2, \tau_1 + \tau_2)$$
in the previous expression, where
$$ \Sigma := \{ (\xi, \tau): (\xi,\tau) = (\xi_1 + \xi_2, \tau_1 + \tau_2) \hbox{ for some } (\xi_1, \tau_1), (\xi_2, \tau_2) \in \Omega \hbox{ such that } |\xi_2| \sim |\xi_1 + \xi_2| \}.$$

We expand this out as
\bas
\int_* &
\chi_{\Omega}(\xi_1, \tau_1)
\chi_{\Omega}(\xi_2, \tau_2)
\chi_\Sigma(\xi_1 + \xi_2, \tau_1 + \tau_2)
\chi_{\Omega}(\xi_3, \tau_3)
\chi_{|\xi_4| \lesssim N_4}
\chi_{|\xi_5| \lesssim N_5}
\chi_{\lambda_4 \lesssim L_4}
\\
&\prod_{j=1}^5 \hat u_j(\xi_j,\tau_j)
\lesssim (N_4 N_5 L_4)^{\frac{1}{2}-\frac{1}{100}} \prod_{j=1}^5 \| u_j \|_{L^2_{x,t}}
\end{align*}
where $\int_*$ denotes integration over the region $\xi_1 + \ldots + \xi_5 = \tau_1 + \ldots + \tau_5 = 0$.  

We can write the left-hand side as\footnote{Alternatively, one could invoke
the Composition Lemma (\cite{tao:xsb}, Lemma 3.7), but the relabeling of indices involved becomes very confusing, so we have chosen to do things explicitly instead.}
\be{ugly}
\int_{**}
\chi_\Sigma(\xi_{12}, \tau_{12})
\chi_{\Omega}(\xi_3, \tau_3)
\chi_{|\xi_4| \lesssim N_4}
\chi_{|\xi_5| \lesssim N_5}
\chi_{\lambda_4 \lesssim L_4}
\hat F(\xi_{12}, \tau_{12}) \prod_{j=3,4,5} \hat u_j(\xi_j,\tau_j) 
\end{equation}
where $\int_{**}$ integrates over the variables $\xi_{12}, \xi_3, \xi_4, \xi_5, \tau_{12}, \tau_3, \tau_4, \tau_5$ such that
$$ \xi_{12} + \xi_3 + \xi_4 + \xi_5 = \tau_{12} + \tau_3 + \tau_4 + \tau_5 = 0,$$
and
$$ \hat F(\xi_{12}, \tau_{12}) := \int_{\xi_1 + \xi_2 = \xi_{12}} \int_{\tau_1 + \tau_2 = \tau_{12}} \chi_{\Omega}(\xi_1,\tau_1) \chi_{\Omega}(\xi_2,\tau_2)
\hat u_1(\xi_1, \tau_1) \hat u_2(\xi_2, \tau_2).$$
Observe that
$$ F = (P_\Omega u_1) (P_\Omega u_2)$$
where $P_\Omega$ is the space Fourier projection corresponding to the set $\Omega$.  From H\"older and \eqref{strich-4} we thus have
$$ \| F \|_{L^2_{x,t}} \leq \| P_\Omega u_1 \|_{L^4_{x,t}} \| P_\Omega u_2 \|_{L^4_{x,t}} \lesssim \| P_\Omega u_1 \|_{0,\frac{1}{3}}
\| P_\Omega u_2 \|_{0,\frac{1}{3}} \lesssim \| u_1 \|_{L^2_{x,t}} \| u_2 \|_{L^2_{x,t}}.$$
It thus suffices to show that the quantity \eqref{ugly} is bounded by
$$ \lesssim (N_4 N_5 L_4)^{\frac{1}{2}-\frac{1}{100}} \| F \|_{L^2_{x,t}} \prod_{j=3,4,5} \| u_j \|_{L^2_{x,t}}.$$

If we relabel $\xi_1, \xi_{23}, \xi_4, \xi_5$ as $\xi'_1, \xi'_2, \xi'_3, \xi'_4$ and similarly for the $\tau$, we thus see that it suffices to show that
$$
\|
\chi_{\Sigma}(\xi'_1, \tau'_1)
\chi_{\Omega}(\xi'_2, \tau'_2)
\chi_{|\xi'_3| \lesssim N'_3}
\chi_{|\xi'_4| \lesssim N'_4}
\chi_{\lambda'_3 \lesssim L'_3}
\|_{[4; \Z \times \R]} \lesssim (N'_3 N'_4 L'_3)^{\frac{1}{2}-\frac{1}{100}}
$$
where we have renamed $N_4, N_5, L_4$ as $N'_3, N'_4, L'_3$ to reduce confusion, and $\lambda'_3$ is short-hand for $\langle \xi'_3 - (\tau'_3)^3 \rangle$.  (Of course we will define the $[4; \Z \times \R]$ norm here to use the primed variables $\xi'_j, \tau'_j$ instead of the unprimed variables.)

Since $|\xi'_3 + \xi'_4| \lesssim N'_3 + N'_4$, we see that the variables $\xi'_1$, $\xi'_2$ are constrained by the relationship $\xi'_1 = -\xi'_2 + O(N'_3 + N'_4)$.  By Schur's test (\cite{tao:xsb}, Lemma 3.11) it thus suffices to show that
\be{bigshot}
\begin{split}
\| \chi_{\Sigma}(\xi'_1,\tau'_1) & \chi_{\Omega}(\xi'_2,\tau'_2)
\chi_{\xi'_1 = -A + O(N'_3 + N'_4)} \chi_{\xi_2 = A + O(N'_3 + N'_4)}\\
&\chi_{|\xi'_3| \lesssim N'_3}
\chi_{|\xi'_4| \lesssim N'_4}
\chi_{\lambda'_3 \lesssim L'_3}
\|_{[4; \Z \times \R]}
\lesssim (N'_3 N'_4 L'_3)^{\frac{1}{2}-\frac{1}{100}}
\end{split}
\end{equation}
for all $A \in \Z$.  

Fix $A$. We may assume that $|A| \gg (N'_3 N'_4 L'_3)^{10}$ since the above expression vanishes otherwise.  

We split into two cases: $L'_3 \leq (N'_3+N'_4)^3$ and $L'_3 \geq (N'_3+N'_4)^3$.

{\bf{Case 1: {$L'_3 \leq (N'_3 + N'_4)^3$ ($L'_3$ not dominant).}}} 

We shall drop the primes from the variables $\xi'_j$, $\tau'_j$.
By Cauchy-Schwarz (\cite{tao:xsb}, Lemma 3.9 or Lemma 3.14) we have
$$ 
\| \chi_{|\xi_2| \lesssim N'_3} \chi_{|\xi_3| \lesssim N'_4}
\chi_{\lambda_2 \lesssim L'_3} \|_{[3; \Z \times \R]}
\lesssim \min(N'_3,N'_4)^{\half} (L'_3)^{\half}$$
so by the Composition Lemma (\cite{tao:xsb}, Lemma 3.7; alternatively one can introduce an ``$F$'' as in the previous arguments) it suffices to show that
$$
\|
\chi_\Sigma(\xi_1,\tau_1)
\chi_\Omega(\xi_2,\tau_2)
\chi_{\xi_1 = -A + O(N'_3 + N'_4)} \chi_{\xi_2 = A + O(N'_3 + N'_4)}
\|_{[3; \Z \times \R]} \lesssim 1$$
since by hypothesis we have
$$ \min(N'_3,N'_4)^{\half} (L'_3)^{\half} 
\lesssim (N'_3 N'_4 L'_3)^{\frac{1}{2}-\frac{1}{100}}.$$
By Cauchy-Schwarz (\cite{tao:xsb}, Lemma 3.9) it suffices to show that
$$
\int \chi_{\Sigma}(\xi_1,\tau_1) \chi_{\Omega}(\xi_2,\tau_2)
\chi_{\xi_1 = -A + O(N'_3 + N'_4)} \chi_{\xi_2 = A + O(N'_3 + N'_4)}\ d\xi_1 d\tau_1
\lesssim 1$$
for all $(\xi_3, \tau_3) \in \Z \times \R$, where $\xi_2, \tau_2$ are given by the formulae $\xi_1 + \xi_2 + \xi_3 = \tau_1 + \tau_2 + \tau_3 = 0$.

Fix $\xi_3$, $\tau_3$.  We may assume that $|\xi_3| \lesssim N'_3 + N'_4$ since the integral vanishes otherwise.  Performing the $\tau$ integral and expanding out $\Sigma$, we reduce to showing that
\begin{align*}
\# \{ (\xi'_1, \xi''_1, \xi_2) \in \Z^3: & \xi'_1 + \xi''_1 + \xi_2 + \xi_3 = 0;
|\xi'_1| \sim A; |\xi''_1| \gg (N'_3 + N'_4)^{10}; \\
&
\xi'_1 + \xi''_1 = -A + O(N'_3 + N'_4); \xi_2 := A + O(N'_3 + N'_4);\\
&{\xi'_1}^3 + {\xi''_1}^3 + \xi_2^3 + \tau_3  = O(1) \}
\lesssim 1.
\end{align*}
By \eqref{4-way} we may rewrite this as
\be{card-2}
\begin{split}
\# \{ (\xi'_1, \xi''_1, \xi_2) \in \Z^3: &\xi'_1 + \xi''_1 + \xi_2 + \xi_3 = 0;
|\xi'_1| \sim A; |\xi''_1| \gg (N'_3 + N'_4)^{10}; \\
&\xi'_1 + \xi''_1 = -A + O(N'_3 + N'_4); \xi_2 := A + O(N'_3 + N'_4); \\
&(\xi'_1 + \xi''_1) (\xi'_1 + \xi_2) (\xi''_1 + \xi_2) \in I \}
\lesssim 1
\end{split}
\end{equation}
where
$$ I := \{ l \in \Z: |l + (\tau_3 - \xi_3^3)/3| \lesssim 1 \}.$$
We may assume that $\tau_3 - \xi_3^3 = O(|A|^3)$, since the left-hand side of \eqref{card-2} vanishes otherwise.  Then from Lemma \ref{div} we see that there are at most $O(1)$ possible values of $\xi'_1 + \xi''_1$ which can contribute to \eqref{card-2}.  Once $\xi'_1 + \xi''_1$ is fixed, we see from elementary algebra that there are at most $O(1)$ triples in \eqref{card} (since $I$ has cardinality $O(1)$), and the claim follows.  This concludes the treatment of Case 1.

{\bf{Case 2: $L'_3 \geq (N'_3+N'_4)^3$ ($L'_3$ dominant).}}

We drop the primes from the variables $\xi'_j, \tau'_j$.
In this case the variable $\tau_3$ is constrained to be $O(L'_3)$, and so it suffices to show that
$$
\| \chi_{\Sigma}(\xi_1,\tau_1) \chi_{\Omega}(\xi_2,\tau_2)
\chi_{\xi_1 = -A + O(N'_3 + N'_4)} \chi_{\xi_2 = A + O(N'_3 + N'_4)}
\chi_{|\tau_3| \lesssim L'_3}
\|_{[4; \Z \times \R]}
\lesssim (N'_3 + N'_4)^{\half},
$$
as the right-hand side is clearly less than $(N'_3 N'_4 L'_3)^{\frac{1}{2}-\frac{1}{100}}$.

By Cauchy-Schwarz (\cite{tao:xsb}, Lemma 3.9) it suffices to show that
$$
\int \chi_{\Sigma}(\xi_1,\tau_1) \chi_{\Omega}(\xi_2,\tau_2)
\chi_{\xi_1 = -A + O(N'_3 + N'_4)} \chi_{\xi_2 = A + O(N'_3 + N'_4)}
\chi_{|\tau_3| \lesssim L'_3}
\ d\xi_1 d\tau_1 d\xi_2 d\tau_2
\lesssim (N'_3 + N'_4)$$
for all $(\xi_4, \tau_4) \in \Z \times \R$, where $\xi_3, \tau_3$ are given by the formulae $\xi_1 + \xi_2 + \xi_3 + \xi_4 = \tau_1 + \tau_2 + \tau_3 + \tau_4 = 0$.  

Fix $\xi_4, \tau_4$.  Performing the $\tau$ integrals and expanding out $\Sigma$, we reduce to showing that
\begin{align*}
\# \{ (\xi'_1, \xi''_1, \xi_2, \xi_3) \in \Z^4: & \xi'_1 + \xi''_1 + \xi_2 + \xi_3 + \xi_4 = 0;
\xi'_1 + \xi''_1 = -A + O(N'_3 + N'_4); \\
& \xi_2 := A + O(N'_3 + N'_4);
{\xi'}_1^3 + {\xi''}_1^3 + \xi_2^3 + \tau_4 = O(L'_3) \}
\lesssim N'_3 + N'_4.
\end{align*}
Note that $\xi_3 + \xi_4 = O(N'_3 + N'_4)$.  From this and \eqref{4-way} we see that it suffices to show that
\be{card-4}
\begin{split}
\# \{ (\xi'_1, \xi''_1, \xi_2, \xi_3) \in \Z^4:& \xi'_1 + \xi''_1 + \xi_2 + \xi_3 + \xi_4 = 0;
\xi'_1 + \xi''_1 = -A + O(N'_3 + N'_4); \\
&
(\xi'_1 + \xi''_1) (\xi''_1 + \xi_2) (\xi_2 + \xi'_1)
\in I \}
\lesssim N'_3 + N'_4
\end{split}
\end{equation}
where
$$ I := \{ l \in \Z: |l + \tau_4/3| \lesssim L'_3 \}.$$
From Lemma \ref{div} we see that there are $O(N'_3 + N'_4)$ possible values of $\xi'_1 + \xi''_1$, the claim follows.

This completes the treatment of Case 2, and hence \eqref{est-1}.  The proof of
Theorem \ref{main} is now complete.
\endprf

\section{Proof of Proposition \ref{main-cor}}\label{cor-sec}

We now prove Proposition \ref{main-cor}.
From Theorem \ref{main}, it suffices to show the bilinear estimate
$$ \| \P(\P(u) \partial_x v) \|_{Z^s} \lesssim \| u \|_{s-1,\half} \|v\|_{s,\half}.$$
Since $\partial_x v = \P(\partial_x v)$ and $\| \partial_x v \|_{s-1,\half} \lesssim \| v\|_{s,\half}$, it suffices to prove the more symmetric estimate
$$ \| \P(\P(u_1) \P(u_2)) \|_{Z^s} \lesssim \| u_1 \|_{s-1,\half} \| u_2 \|_{s-1,\half}.$$
The estimate
\be{kpv-est}
\| \P(\P(u_1) \P(u_2)) \|_{s,-\half}
\sim \| (\P(u_1) \P(u_2))_x \|_{s-1,-\half}
\lesssim \| u_1 \|_{s-1,\half} \| u_2 \|_{s-1,\half}
\end{equation}
is proven in \cite{kpv:kdv} (see also \cite{tao:xsb}, Corollary 6.5).  Thus by \eqref{z-def} it remains only to show
$$
\| 
\frac{\langle \xi \rangle^s \chi_{\xi \neq 0} \widehat{\P(u_1) \P(u_2)}(\xi,\tau)}{\langle \tau - \xi^3\rangle} \|_{L^2_\xi L^1_\tau}
\lesssim 
\| u_1 \|_{s-1,\half} \| u_2 \|_{s-1,\half}.
$$
We shall actually prove the stronger estimate
\be{cor-plus}
\| 
\frac{\langle \xi \rangle^s \chi_{\xi \neq 0} \widehat{\P(u_1) \P(u_2)}(\xi,\tau)}{\langle \tau - \xi^3\rangle^{1-\delta}} \|_{L^2_\xi L^1_\tau}
\lesssim 
\| u_1 \|_{s-1,\half} \| u_2 \|_{s-1,\half}
\end{equation}
for some small $0 < \delta \ll 1$.

The difficulty here is that
$$ \int \frac{d\tau}{\langle \tau - \xi^3\rangle^{1-2\delta}}$$
is divergent for each $\xi$, since otherwise we could use Cauchy-Schwarz in $\tau$ to reduce this to \eqref{kpv-est}.  Our strategy shall then be to somehow mollify the weight $\frac{1}{\langle \tau - \xi^3 \rangle^{1-\delta}}$ so that the above integral is no longer divergent.

We may assume that the spacetime Fourier transforms of $u_1$ and $u_2$ are non-negative.  We expand out the left-hand side as
$$
\| 
\sum_{\xi_1} \int_{\tau_1}
\frac{\langle \xi \rangle^s \chi_{\xi_1 \xi_2 \xi \neq 0}
\hat u_1(\xi_1,\tau_1) \hat u_2(\xi_2,\tau_2)}
{\langle \tau - \xi^3\rangle^{1-\delta}}\ d\tau_1 \|_{L^2_\xi L^1_\tau}
$$
where $\xi_1 + \xi_2 = \xi$ and $\tau_1 + \tau_2 = \tau$.

From the estimate
$$ \frac{\langle \xi \rangle^s}{\langle \xi_1 \rangle^{s-1} \langle \xi_2 \rangle^{s-1}} \lesssim |\xi_1|^{\half} |\xi_2|^{\half} |\xi|^{\half}$$
when $\xi_1 \xi_2 \xi \neq 0$ and $\xi = \xi_1 + \xi_2$, we reduce to
\be{lhs}
\| 
\sum_{\xi_1} \int_{\tau_1}
\frac{|\xi_1|^{\half} |\xi_2|^{\half} |\xi_3|^{\half} \hat u_1(\xi_1,\tau_1) \hat u_2(\xi_2,\tau_2)}
{
\langle \tau - \xi^3\rangle^{1-\delta}
\langle \tau_1 - \xi_1^3 \rangle^{\half}
\langle \tau_2 - \xi_2^3 \rangle^{\half}
}\ d\tau_1 \|_{L^2_\xi L^1_\tau}
\lesssim \| u_1 \|_{0,0} \| u_2 \|_{0,0}.
\end{equation}

From the identity
$$
(\tau - \xi^3) = (\tau_1 - \xi_1^3) + (\tau_2 - \xi_2^3) - 3 \xi_1 \xi_2 \xi$$
we see that at least one of the quantities
$$ \langle \tau - \xi^3 \rangle, \langle \tau_1 - \xi_1^3 \rangle, \langle \tau_2 - \xi_2^3 \rangle$$
must exceed $\gtrsim |\xi_1| |\xi_2| |\xi|$.
Suppose that we had
$$ \langle \tau_1 - \xi_1^3 \rangle \gtrsim |\xi_1| |\xi_2| |\xi|.$$
Then we can reduce \eqref{lhs} to
$$
\| 
\langle \tau - \xi^3  \rangle^{-1+\delta}
\sum_{\xi_1} \int_{\tau_1}
\hat u_1(\xi_1,\tau_1) \hat u_2(\xi_2,\tau_2)\ d\tau_1 \|_{L^2_\xi L^1_\tau}
\lesssim \| u_1 \|_{0,0} \| u_2 \|_{0,\half}.
$$
If $\delta$ is sufficiently small, the weight $\langle \tau - \xi^3 
\rangle^{-4/3+2\delta}$ is integrable in $\tau$ uniformly in $\xi$.  By H\"older in $\tau$ it thus suffices to show
$$
\| 
\langle \tau - \xi^3  \rangle^{-\frac{1}{3}}
\sum_{\xi_1} \int_{\tau_1}
\hat u_1(\xi_1,\tau_1) \hat u_2(\xi_2,\tau_2)\ d\tau_1 \|_{L^2_\xi L^2_\tau}
\lesssim \| u_1 \|_{0,0} \| u_2 \|_{0,\half},
$$
or equivalently that
\be{red}
|\int u_1 u_2 u_3\ dx dt| \lesssim \| u_1 \|_{0,0} \| u_2 \|_{0,\half} \| u_3 \|_{0,\frac{1}{3}}.
\end{equation}
But this follows from H\"older and \eqref{strich-4}.
From the above and symmetry, we may reduce to the case where
$$ \langle \tau - \xi^3  \rangle \gtrsim |\xi_1| |\xi_2| |\xi|.$$
Suppose for the moment that we also had
$$ \langle \tau_1 - \xi_1^3  \rangle \gtrsim (|\xi_1| |\xi_2| |\xi|)^{\frac{1}{100}}.$$
Then we can reduce \eqref{lhs} to
$$
\| 
\langle \tau - \xi^3  \rangle^{-\half-\frac{1}{600}+\delta}
\sum_{\xi_1} \int_{\tau_1}
\hat u_1(\xi_1,\tau_1) \hat u_2(\xi_2,\tau_2)\ d\tau_1 \|_{L^2_\xi L^1_\tau}
\lesssim \| u_1 \|_{0,\frac{1}{3}} \| u_2 \|_{0,\half}.
$$
If $\delta$ is sufficiently small, $\langle \tau - \xi^3  \rangle^{-1-\frac{1}{300}+2\delta}$ is integrable in $\tau$ uniformly in $\xi$, and we may use H\"older in $\tau$ and duality to reduce to \eqref{red} as before.

We may thus assume that
$$ \langle \tau_i - \xi_i^3  \rangle \ll (|\xi_1| |\xi_2| |\xi|)^{\frac{1}{100}}$$
for $i=1,2$.  

In particular, we have that $\xi_1, \xi_2, \xi \neq 0$ and
$$
\tau - \xi^3  = -3 \xi_1 \xi_2 \xi + O(\langle \xi_1 \xi_2 \xi\rangle)^{
\frac{1}{100}}
$$
and hence that
$$ \langle \tau - \xi^3 \rangle \sim \langle \xi_1 \xi_2 \xi \rangle.$$
Applying these estimates and changing the $\tau_1, \tau$ integrals to $\tau_1, \tau_2$ integrals, we may therefore majorize the left-hand side of \eqref{lhs} by
$$
\| 
\sum_{\xi_1} \langle \xi_1 \xi_2 \xi \rangle^{\delta - \half}
\int_{\tau_1 = \xi_1^3 + O(\langle \xi_1 \xi_2 \xi \rangle^{\frac{1}{100}})} \int_{\tau_2 = \xi_2^3 + O(\langle \xi_1 \xi_2 \xi \rangle^{\frac{1}{100}})}
\hat u_1(\xi_1,\tau_1) \hat u_2(\xi_2,\tau_2)\ d\tau_1 d\tau_2 \|_{L^2_\xi}.
$$
Applying Cauchy-Schwarz in $\tau_1$ and $\tau_2$ separately, we may majorize this by
$$\| 
\sum_{\xi_1} \langle \xi_1 \xi_2 \xi \rangle^{\delta + \frac{1}{100} - \half}
F_1(\xi_1) F_2(\xi_2) \|_{L^2_\xi}$$
where
$$ F_i(\xi) := (\int \hat u_i(\xi,\tau)^2\ d\tau)^{\half}.$$
Since $\xi_1 + \xi_2 = \xi$, we have
$$ \langle \xi_1 \xi_2 \xi \rangle^{\delta + \frac{1}{100} - \half} \lesssim \langle \xi \rangle^{-\half-}$$
if $\delta$ is sufficiently small.  Thus by H\"older we may majorize the previous by
$$ \| \sum_{\xi_1} F_1(\xi_1) F_2(\xi_2) \|_{L^\infty_\xi}.$$
But by Cauchy-Schwarz this is bounded by $\|F_1\|_2 \|F_2\|_2 = \|u_1\|_{0,0} \|u_2\|_{0,0}$ as desired.
\endprf

\section{Local well-posedness for small data}\label{local-sec}

We can now prove Theorem \ref{periodic-kdv}.  Let $s \geq \half$, and let $u_0$ be initial data with small $H^s$ norm.  Write \eqref{pkdv} as
$$
u_t + \frac{1}{4\pi^2} u_{xxx} + F'(u) u_x = 0; \quad u(x,0) = u_0(x)
$$
where $F'$ is the derivative of the polynomial $F$.

We now follow the standard reductions of \cite{staff:pkdv}.  We apply the ``gauge transformation''
$$ v(x,t) := u(x - (\int_0^t \int_\T F'(u)(x',t')\ dx' dt'), t).$$
This transformation is invertible:
$$ u(x,t) := v(x + (\int_0^t \int_\T F'(v)(x',t')\ dx' dt'), t).$$
Also, it preserves the initial data $u_0$, and is a homeomorphism on $H^s(\T)$ for each time $t$.  (Note from Sobolev embedding and the hypothesis $s \geq \half$ that $F(u)$ is locally integrable whenever $u \in H^s(\T)$).  

It is easy to check that $u$ solves \eqref{pkdv} if and only if $v$ solves the equation
$$
v_t + \frac{1}{4\pi^2} v_{xxx} + \P(F'(v)) v_x = 0; \quad v(x,0) = u_0(x).$$
Since $F'(v) v_x = F(v)_x$ and $v_x$ both have mean zero, $\P(F'(v)) v_x$ must also have mean zero.  Thus we may rewrite the above Cauchy problem as
$$
v_t + \frac{1}{4\pi^2} v_{xxx} + \P(\P(F'(v)) v_x) = 0; \quad v(x,0) = u_0(x),
$$
or in integral form as
$$
v(t) = S(t) u_0 - \int_0^t S(t-t') \P(\P(F'(v)) v_x)(t')\ dt';
$$
recall that
$$S(t) := \exp(-\frac{1}{4\pi^2} t \partial_{xxx})$$
is the fundamental solution of the Airy equation.  If we are only interested in solving this equation up to time 1, we may (following \cite{borg:xsb}) replace this equation\footnote{Alternatively, one can restrict time to $[0,1]$ and replace the $X^{s,b}$ norms by their equivalence class counterparts on this time interval.} with
$$
v(t) = \eta(t) (S(t) u_0 - \int_0^t S(t-t') \P(\P(F'(v)) v_x)(t')\ dt')
$$
where $\eta$ is the bump function from Section \ref{embedding-sec}. 

We shall apply the contraction mapping principle to the map
\be{contraction}
v \mapsto \eta(t) (S(t) u_0 - \int_0^t S(t-t') \P(\P(F'(v)) v_x)(t')\ dt').
\end{equation}

From several applications of Proposition \ref{main-cor} we have
$$ \| \P(\P(F'(v)) v_x) - \P(\P(F'(w)) w_x) \|_{Z^s} 
\ll \| v - w \|_{Y^s}$$
if $\|v\|_{Y^s}, \|w\|_{Y^s}$ are sufficiently small.  Also, it is easily verified that
$$ \| \eta(t) F \|_{Z^s} \lesssim \| F\|_{Z^s}$$
for all $F$.  From these estimates and \eqref{grow}, Lemma \ref{invert} we see that \eqref{contraction} is a contraction on a small ball of $Y^s$
if $\| u_0 \|_{H^s}$ is sufficiently small.  This will gives local existence, continuity, and uniqueness in the space $Y^s$, which embeds into $C([0,1]; H^s)$ by \eqref{sup-est}.  The proof of Theorem \ref{main} is now complete.
\endprf

\begin{remark}
\label{largedata}
Note that this argument can be modified to deal with the large data case.  The point is that in many of the above estimates, at least one of the $X^{s,\half}$ norms can be replaced with a $X^{s,\half-\delta}$ norm (the exact choice of factor may depend on what case one is in).  If one localizes to a small time $T$, one can estimate the $X^{s,\half-\delta}$ norm by the $X^{s,\half}$ norm and gain a small power of $T$.  This allows one to obtain the desired contraction if $T$ is sufficiently small depending on the $H^s$ norm of the initial data.  See e.g. \cite{selberg:thesis}, or equation (3.9) in \cite{cst}.  However, we shall not pursue these arguments, and rely instead on the rescaling arguments of the following sections (which automatically give the correct power dependence of $T$ on $\|u_0\|_{H^s}$; this seems quite difficult to do using the $X^{s,\half-\delta}$ norms if one refuses to rescale).
\end{remark}

\section{Large periods}\label{period-sec}

We now begin the proof of Theorem \ref{gwp}.  The first step is to use rescaling arguments to generalize the previous estimates to the large period case.  In order to do this we shall need to set up some conventions for Fourier transforms, $X^{s,b}$ spaces, etc. in the large period case.  (These conventions are also used in \cite{ckstt:2}).

Fix $\lambda \gg 1$.  

In the sequel $(d\xi)_\lambda$ will be normalized counting measure on $\Z/\lambda$:
$$ \int a(\xi) (d\xi)_\lambda := \frac{1}{\lambda} \sum_{\xi \in \Z/\lambda} a(\xi).$$
Thus $(d\xi)_\lambda$ is the counting measure on the integers when $\lambda = 1$, and converges weakly to Lebesgue measure when $\lambda \to \infty$.  

In the remainder of this section, all Lebesgue norms in $\xi$ will be with respect to the measure $(d\xi)_\lambda$, while all Lebesgue norms in $x$ will be on the large torus $\R/\lambda \Z$.

Let $u(x,t)$ be a function of $\R/\lambda \Z \times \R$.  We define the spacetime Fourier transform $\hat u$ to be the function
$$
\hat u(\xi,\tau) := \int_{\R/\lambda \Z} e^{-2\pi i (\xi x + t\tau)} f(x)\ dx
$$
defined for all $\xi \in \Z/\lambda$.  The inverse Fourier transform is given by
$$
u(x,t) = \int \hat u(\xi,t) e^{2\pi i (\xi x + t\tau)} (d\xi)_\lambda d\tau.
$$
We define the spatial Fourier transform $\hat f(\xi)$ similarly.

We define the Sobolev spaces $H^s_\lambda$ on $[0,\lambda]$ by
$$
\| f \|_{H^s_\lambda} := \| \hat f(\xi) \langle \xi \rangle^s \|_{L^2_\xi}$$
and the spaces $X^{s,b}_\lambda$ on $[0, \lambda] \times \R$ by
$$
\| u \|_{X^{s,b}_\lambda} := \| \hat u(\xi, \tau) \langle \xi \rangle^s \langle \tau - \xi^3  \rangle^b \|_{L^2_{\tau,\xi}}.$$
We also define the spaces $Y^s_\lambda$, $Z^s_\lambda$ as 
$$ \| u \|_{Y^s_\lambda} := \| u \|_{X^{s,\half}_\lambda} + \| \langle \xi \rangle_\lambda^s \hat u \|_{L^2_\xi L^1_\tau}$$
$$ \| u \|_{Z^s_\lambda} := \| u \|_{X^{s,-\half}_\lambda} + \| 
\frac{\langle \xi \rangle_\lambda^s \hat u}{\langle \tau - \xi^3  \rangle} \|_{L^2_\xi L^1_\tau}.$$

\begin{remark}
\label{largedatatwo}
Our strategy for the large data theory will be to rescale large $H^s$ data in the period 1 case to small $H^s_\lambda$ data in the period $\lambda$ case, for some large $\lambda$ depending on the norm of the original data.  This procedure works well when $k=3$, but runs into a difficulty when $k \geq 4$ since the $L^2$ component of the $H^s_\lambda$ norm is critical or supercritical.  This difficulty can probably be avoided by modifying the Fourier weight of the $H^s_\lambda$,  $X^{s,b}_\lambda$, $Y^s_\lambda$, $Z^s_\lambda$ spaces at low frequencies, but we shall not discuss these matters here, and focus instead on the $k=3$ case. A related problem has been addressed \cite{ckstt:1} in the context of cubic NLS on $\R^3$ where the $L^2$ norm is supercritical.
\end{remark}

Not all of the embeddings in \eqref{embedding-sec} still hold.  However, we have the analogue
\be{strich-lambda}
\| u \|_{L^4_{x,t}} \lesssim \| u \|_{X^{0,\frac{1}{3}}_\lambda}
\end{equation}
of \eqref{strich-4}, which just follows from rescaling\footnote{Observe that if we give $x$, $t$ the units of $length$ and $length^3$ respectively, so that $\xi$ and $\tau$ have units $length^{-1}$ and $length^{-3}$, then both sides have the units of $length$ and thus scale properly (the fact that we have the weight $\langle \tau - \xi^3 \rangle$ instead of $|\tau - \xi^3|$ affects this slightly, but the effect is favorable).}  \eqref{strich-4}.  Also, we have 
\be{infty-lambda}
\| u \|_{L^\infty_{x,t}} \lesssim \| u \|_{X^{\half+,\half+}_\lambda}
\end{equation}
\be{sobolev-lambda}
\| u \|_{L^2_t L^\infty_x} \lesssim \| u \|_{X^{\half+,0}_\lambda}
\end{equation}
\be{energy-lambda}
\| u \|_{L^\infty_t L^2_x} \lesssim \| u \|_{X^{0,\half+}_\lambda};
\end{equation}
These estimates are proved in exactly the same way as their $\lambda = 1$ counterparts.

We now develop analogues of the preceding results for large $\lambda$.  We begin with the analogue of Theorem \ref{main}.

\begin{proposition}\label{main-lambda}
For all $\half \leq s \leq 1$ and $k \geq 3$, we have
$$ \| u_1 \ldots u_k \|_{X^{s-1,\half}_\lambda} \lesssim 
\lambda^{0+} 
\prod_{i=1}^k \| u_i \|_{Y^s_\lambda}.$$
\end{proposition}

\begin{proof}  
We repeat the reductions in Section \ref{Non-endpoint-sec}.  We first consider the contribution of the case \eqref{tau1-dom}.  In this case it suffices to show
$$ \| u_1 \ldots u_k \|_{X^{s-1,0}_\lambda} \lesssim 
\| u_1 \|_{X^{s,0}_\lambda} 
\prod_{i=2}^k \| u_2 \|_{X^{s,\half}_\lambda}.$$
On the other hand, from the Sobolev embeddings $H^s_\lambda \subseteq L^{2k}$, $L^2 \subseteq H^{s-1}_\lambda$ and H\"older we have the spatial estimate
$$ \| f_1 f_2 f_3 \|_{H^{s-1}_\lambda} \lesssim 
\prod_{i=1}^k \| f_i \|_{H^s_\lambda},
$$
and the claim follows by setting $f_i = u_i(t)$ and then taking $L^2$ norms in time.

By symmetry it remains only to consider the case \eqref{tau-dom}. By the arguments of Section \ref{Non-endpoint-sec} we thus reduce to
\be{est-2-lambda}
\| \frac{|\xi_1|^{\half} |\xi_2|^{\half} |\xi_{k+1}|^{\half}}{
(\prod_{i=1}^k \langle \xi_i \rangle^s \langle \tau_i - \xi_i^3 \rangle^{\half})
\langle \xi_{k+1} \rangle^{1-s}
}
\|_{[4; \Z/\lambda \times \R]} \lesssim \lambda^{0+}
\end{equation}
and 
\be{est-1-lambda}
\| \frac{|\xi_1|^{\half} |\xi_2|^{\half} |\xi_3|^{\half}}{
(\prod_{i=1}^k \langle \xi_i \rangle^s \langle \tau_i - \xi_i^3 \rangle^{\half})
\langle \xi_{k+1} \rangle^{1-s}
}
\|_{[4; \Z/\lambda \times \R]} \lesssim \lambda^{0+}.
\end{equation}

We may assume that $|\xi_1| \geq \ldots \geq |\xi_k|$ for these estimates.  This implies that $|\xi_4| \lesssim |\xi_1|$.  In particular, the $s > \half$ form of these estimates will then follow from the $s=\half$ case and the Comparison Principle (\cite{tao:xsb}, Lemma 3.1).  We shall thus assume $s=\half$ in the sequel.

Consider \eqref{est-2-lambda}.  As in Section \ref{2-prf}, it suffices to show
\be{e2l}
\| u_1 \ldots u_k \|_{L^2_{x,t}} \lesssim \lambda^{0+}
\| u_1\|_{X^{0,\half}_\lambda} \| u_2 \|_{X^{0,\half}_\lambda} 
\prod_{i=3}^k \| u_i \|_{X^{\half, \half}_\lambda}.
\end{equation}

By rescaling \eqref{est2-plus} and conceding several powers of $\lambda$ we obtain
$$ \| u_1 \ldots u_k \|_{L^2_{x,t}} \lesssim \lambda^C
\| u_1\|_{X^{0,\half-\delta}_\lambda} \| u_2 \|_{X^{0,\half-\delta}_\lambda} 
\prod_{i=3}^k \| u_i \|_{X^{\half-\delta,\half-\delta}_\lambda}$$
for some large constant $C$.  On the other hand, from \eqref{strich-lambda}, \eqref{infty-lambda} we have
$$ \| u_1 \ldots u_k \|_{L^2_{x,t}} \lesssim 
\| u_1\|_{X^{0,\frac{1}{3}}_\lambda} \| u_2 \|_{X^{0,\frac{1}{3}}_\lambda} 
\prod_{i=3}^k \| u_i \|_{X^{\half+,\half+}_\lambda}.$$

If one interpolates this estimate a little bit with the previous one, one obtains \eqref{e2l} as desired (in fact we even get a gain in some of the  indices).

Now consider \eqref{est-1-lambda}.  As in Section \ref{1-prf}, it suffices to show
$$|\int u_1 \ldots u_{k+1}\ dx dt| \lesssim \lambda^{0+}
\| u_1 \|_{X^{0,\half}_\lambda}
\| u_2 \|_{X^{0,\half}_\lambda }
\| u_3 \|_{X^{0,\half}_\lambda }
(\prod_{i=4}^k \| u_i \|_{X^{\half,\half}_\lambda})
\| u_{k+1} \|_{X^{\half,0}_\lambda }
$$
Arguing as with \eqref{est-2-lambda}, using \eqref{est1-plus} instead of \eqref{est2-plus}, we reduce to showing that
$$|\int u_1 \ldots u_{k+1}\ dx dt| \lesssim 
\| u_1 \|_{X^{0,\half+}_\lambda}
\| u_2 \|_{X^{0,\half+}_\lambda }
\| u_3 \|_{X^{0,\half+}_\lambda }
(\prod_{i=4}^k \| u_i \|_{X^{\half +,\half +}_\lambda})
\| u_{k+1} \|_{X^{\half +,0}_\lambda }.$$
But this follows from two applications of \eqref{strich-lambda}, one application each of \eqref{energy-lambda} and \eqref{sobolev-lambda}, and $k-3$ applications of \eqref{infty-lambda}.
\end{proof}

In \cite{ckstt:2} the following large-period analogue of \eqref{kpv-est} was proven:
\begin{proposition}\label{kpv-large}(\cite{ckstt:2}, equation (7.34))
We have
$$
|\int \P(u_1) \P(u_2) \P(u_3)\ dx dt| \lesssim
\lambda^{0+}
\| u_1 \|_{-\half,\half}
\| u_2 \|_{-\half,\half}
\| u_3 \|_{-\half,\half}
$$
\end{proposition}

Combining this estimate with Proposition \ref{main-lambda}
we obtain
\begin{corollary}\label{endpoint}
We have the quintilinear estimate 
$$ | \int \P(u_1 u_2 u_3) \P(u_4) \P(u_5)\ dx dt|
\lesssim \lambda^{0+}
\| u_1 \|_{Y^{\half}_\lambda}
\| u_2 \|_{Y^{\half}_\lambda}
\| u_3 \|_{Y^{\half}_\lambda}
\| u_4 \|_{Y^{-\half}_\lambda}
\| u_5 \|_{Y^{-\half}_\lambda}.
$$
\end{corollary}

This estimate is required in \cite{ckstt:2} to prove global well-posedness of the KdV and modified KdV equation for $s \geq -\half $ and $s \geq \half$ respectively.

From Propositions \ref{main-lambda} and \ref{kpv-large} we may also deduce a $k=3$ rescaled version of Corollary \ref{main-cor}.
\begin{corollary}\label{rescaled}
We have
\be{lwp}
\| \P(\P(u_1 u_2 u_3) \partial_x u_4) \|_{Z^s_\lambda} \lesssim \lambda^{0+} 
\| u_1 \|_{Y^s_\lambda}
\| u_2 \|_{Y^s_\lambda}
\| u_3 \|_{Y^s_\lambda}
\| u_4 \|_{Y^s_\lambda},
\end{equation}
for all $\half \leq s \leq 1$.
\end{corollary}

\begin{proof}
By Proposition \ref{main-lambda} and the observation that $\partial_x u_4 = \P(\partial_x u_4)$ obeys the estimate
$$ \| \partial_x u_4\|_{Y^{s-1}_\lambda} \lesssim \| u_4 \|_{Y^s_\lambda} $$
it suffices to show that
$$
\| \P( \P(u) \P(v) ) \|_{Z^s_\lambda} \lesssim \lambda^{0+}
\| u \|_{X^{s-1,\half}_\lambda} \| v\|_{Y^{s-1}_\lambda}.$$

The $X^{s,-\half}_\lambda$ portion of $Z^s_\lambda$ is acceptable by Proposition \ref{kpv-large} and duality.  It thus suffices to show that
\be{cpl}
\| \langle \xi \rangle^s  \chi_{\xi \neq 0}
\frac{\widehat{\P(u) \P(v)}(\xi,\tau)}{\langle \xi - \tau^3 \rangle}
\|_{L^2_\xi L^1_\tau} \lesssim \lambda^{0+}
\| u \|_{X^{s-1,\half}_\lambda} \| v\|_{Y^{s-1}_\lambda}.
\end{equation}
By rescaling \eqref{cor-plus} we see that
$$
\| \langle \xi \rangle^s  \chi_{\xi \neq 0}
\frac{\widehat{\P(u) \P(v)}(\xi,\tau)}{\langle \xi - \tau^3 \rangle^{1-\delta}}
\|_{L^2_\xi L^1_\tau} \lesssim \lambda^{C}
\| u \|_{X^{s-1,\half}_\lambda} \| v\|_{Y^{s-1}_\lambda}$$
for some large constant $C$.  It will thus suffice to show that
$$
\| \langle \xi \rangle^s  \chi_{\xi \neq 0}
\frac{\widehat{\P(u) \P(v)}(\xi,\tau)}{\langle \xi - \tau^3 \rangle^{1+}}
\|_{L^2_\xi L^1_\tau} \lesssim \lambda^{0+}
\| u \|_{X^{s-1,\half}_\lambda} \| v\|_{Y^{s-1}_\lambda},$$
since \eqref{cpl} then follows by interpolating this estimate a little bit with the previous one.  But by Cauchy-Schwarz in $\tau$ we may estimate the left-hand side of the above estimate by
$$
\| \langle \xi \rangle^s  \chi_{\xi \neq 0}
\frac{\widehat{\P(u) \P(v)}(\xi,\tau)}{\langle \xi - \tau^3 \rangle^{\half}}
\|_{L^2_\xi L^2_\tau},$$
and the claim follows as before from Proposition \ref{kpv-large} and duality.
\end{proof}

\section{An interpolation lemma}\label{interp-sec}

The purpose of this section is to prove a general interpolation result which will be useful in low regularity global well-posedness theory.  Roughly speaking, this result asserts that if one can prove local well-posedness at two different levels of regularity $H^{s_1}$ and $H^{s_2}$ with $s_1 > s_2$, then one can also prove local well-posedness at a regularity which behaves like $H^{s_1}$ for low frequencies and $H^{s_2}$ for high frequencies.

We need some notation.  Let $m(\xi)$ be a smooth non-negative symbol on $\R$ which equals $1$ for $|\xi| \leq 1$ and equals $|\xi|^{-1}$ for $|\xi| \geq 2$.

For any $N \geq 1$ and $\alpha \in \R$, let $I_N^\alpha$ denote the spatial Fourier multiplier
$$ \widehat{I_N^\alpha f}(\xi) = m(\frac{\xi}{N})^\alpha \hat f(\xi).$$
The operator $I_1^\alpha$ is thus a standard smoothing operator of order $\alpha$.  The operators $I_N^\alpha$ are similar operators, but are somewhat larger; for instance, $I_N^\alpha$ is the identity on low frequencies $|\xi| \leq N$.

One can of course apply these operators to spacetime functions $u(x,t)$ by the formula
$$ \widehat{I_N^\alpha u}(\xi,\tau) = m(\frac{\xi}{N} )^\alpha \hat u(\xi,\tau).$$

\begin{definition}\label{additive-def}  For every $x$, let $S_x$ denote the shift operator $S_x u(x',t) := u(x'-x,t)$.  A Banach space $X$ of spacetime functions is said to be \emph{translation invariant} if one has 
$$ \| S_x u \|_X = \| u \|_X$$
for all $x$ and all $u \in X$.

A multilinear operator $T(u_1, \ldots, u_n)$ is said to be \emph{translation invariant} if one has
$$ S_x T(u_1, \ldots, u_n) = T(S_x u_1, \ldots, S_x u_n)$$
for all $x$.
\end{definition}

Equivalently, an operator $T$ is translation invariant if its kernel $K(x,t; x_1,t_1, \ldots, x_n,t_n)$ (in the sense of distributions) has the symmetry 
$$ K(x,t; x_1,t_1, \ldots, x_n,t_n) = K(x+y, t; x_1 + y, t_1, \ldots, x_n + y, t_n).$$
From this one can easily see that if $T$ is translation invariant and each $u_i$ has Fourier support in the region $\{ (\xi_i, \tau_i): \xi_i \in \Omega_i \}$ for some sets $\Omega_i$, then $T(u_1, \ldots, u_n)$ must have Fourier support in the Minkowski sum $\{ (\xi, \tau): \xi \in \Omega_1 + \ldots + \Omega_n \}$ of the above regions.

If a Banach space $X$ is translation invariant, then $X$ is closed under convolutions with $L^1$ kernels.  In particular, Littlewood-Paley projection operators are bounded on $X$.

\begin{lemma}\label{interp-lemma}
Let $\alpha_0 > 0$ and $n \geq 1$.  Suppose that $Z$, $X_1, \ldots, X_n$ are translation invariant Banach spaces and $T$ is a translation invariant $n$-linear operator such that has the estimate
\be{interp-hyp}
\| I_1^\alpha T(u_1, \ldots, u_n) \|_Z \lesssim \prod_{i=1}^n \| I_1^\alpha u_i \|_{X_i}
\end{equation}
for all $u_1, \ldots, u_n$ and all $0 \leq \alpha \leq \alpha_0$.  Then one has the estimate
\be{interp-conc}
\| I_N^\alpha T(u_1, \ldots, u_n) \|_Z \lesssim \prod_{i=1}^n \| I_N^\alpha u_i \|_{X_i}
\end{equation}
for all $u_1, \ldots, u_n$, all $0 \leq \alpha \leq \alpha_0$, and $N \geq 1$, with the implicit constant independent of $N$.
\end{lemma}

Note that \eqref{interp-conc} trivially follows from \eqref{interp-hyp} when $N=1$ or $N=\infty$.  We remark for future applications that this lemma works in general spatial dimension, and in both the periodic and non-periodic cases (with any period $\lambda \geq 1$).

\begin{proof}
Since Littlewood-Paley  projection operators are bounded on the $X_i$, we may split each $u_i$ into a piece supported on frequencies $|\xi_i| \lesssim N$ and a piece supported on frequencies $|\xi_i| \gg N$, and deal with each contribution separately.

First suppose that each $u_i$ has Fourier support on $|\xi_i| \lesssim N$.  Since $T$ is translation invariant, $T(u_1, \ldots, u_n)$ also has Fourier support on a region $|\xi| \lesssim N$.  On these regions the operators $I_N^{\pm \alpha}$ are essentially Littlewood-Paley multipliers, and we thus have 
$$\| I_N^\alpha T(u_1, \ldots, u_n) \|_Z \lesssim \| T(u_1, \ldots, u_n) \|_Z$$
and
$$ \| u_i \|_{X_i} \lesssim \| I_N^{\alpha} u_i\|_{X_i}.$$
The claim then follows from the $\alpha = 0$ version of \eqref{interp-hyp}.

Now suppose that $u_1$ (for instance) has Fourier support on $|\xi_1| \gg N$.  Then we have
$$ \| I_1^\alpha u_1 \|_{X_1} = N^{-\alpha} \| I_N^\alpha u_1 \|_{X_1}.$$
Also, since the operator $I_1^\alpha I_N^{-\alpha}$ has an integrable kernel, we have
$$ \| I_1^\alpha u_i \|_{X_i} \lesssim \| I_N^\alpha u_i \|_{X_i}$$
for $1 < i \leq n$.  By \eqref{interp-hyp} we thus have
$$ \| I_1^\alpha T(u_1, \ldots, u_n) \|_Z \lesssim N^{-\alpha} \prod_{i=1}^n \| I_N^\alpha u_i \|_{X_i}.$$

Since $N^{-\alpha} I_N^\alpha I_1^{-\alpha}$ has an integrable kernel, the claim then follows.
\end{proof}

From Corollary \ref{rescaled} and Lemma \ref{interp-lemma} (with $\alpha_0 = \half$) we have in particular that

\begin{corollary}\label{rescaled-I}
We have
$$
\| I_N^{1-s} \P(\P(u_1 u_2 u_3) \partial_x u_4) \|_{Z^1_\lambda} \lesssim \lambda^{0+} 
\| I_N^{1-s} u_1 \|_{Y^1_\lambda}
\| I_N^{1-s} u_2 \|_{Y^1_\lambda}
\| I_N^{1-s} u_3 \|_{Y^1_\lambda}
\| I_N^{1-s} u_4 \|_{Y^1_\lambda},
$$
for all $\half \leq s \leq 1$.
\end{corollary}

\section{Global well-posedness}\label{global-sec}

We now give the proof of Theorem \ref{gwp}, following the general ``I-method'' scheme in \cite{keel:mkg}, \cite{ckstt:1}, \cite{ckstt:2}, \cite{ckstt:4} (see also \cite{keel:wavemap}).  

Fix $\half \leq s$.  In light of the results of \cite{staff:pkdv} we may assume that $s < 1$.  Fix $u_0 \in H^s$; the norm $\| u_0 \|_{H^s}$ may possibly be large.  We shall need an absolute constant $0 < \eps \ll 1$ to be chosen later.  We shall also need parameters $\lambda \gg 1$, $N \gg 1$ depending on $\eps$ and $\| u_0 \|_{H^s}$ to be chosen later.

As in Section \ref{local-sec}, it suffices to find an $H^s$ solution $v(x,t)$ to the problem 
$$
v_t + \frac{1}{4\pi^2} v_{xxx} + \P(\P(v^3) v_x) = 0; \quad v(x,0) = u_0(x).
$$

The first step is to rescale the problem.  Consider the functions $u_0^\lambda$, $v^\lambda$ on $\R / \lambda \Z \times \R$ defined by
$$ u_0^\lambda(x) := \lambda^{-\frac{2}{3}} u_0(\frac{x}{\lambda}), \quad
v^\lambda(x,t) := \lambda^{-\frac{2}{3}} v(\frac{x}{\lambda}, \frac{t}{\lambda^3}).$$
We shall construct a solution to the Cauchy problem
\be{rescaled-cauchy}
v^\lambda_t + \frac{1}{4\pi^2} v^\lambda_{xxx} + \P(\P((v^\lambda)^3) v^\lambda_x) = 0; \quad v^\lambda(x,0) = u^\lambda_0(x)
\end{equation}
up to some time $T > 0$ to be determined later; this will yield a solution to the original problem up to time $\frac{T}{\lambda^3}$.

The $L^2$ norm is sub-critical:
$$ \| u_0^\lambda \|_2 = \lambda^{-\frac{1}{6}} \| u_0 \|_2.$$
Thus if $\lambda$ is sufficiently large depending on $\| u_0 \|_{H^s}$ and $\eps$ we have
\be{sev-bound}
\| u_0^\lambda \|_2 \leq \eps.
\end{equation}
For any function $v$ on $\R/\lambda \Z$, define the Hamiltonian $H(v)$ by
$$ H(v) := \int \frac{1}{8 \pi^2} v_x^2 - \frac{1}{20} v^5\ dx.$$
One can verify that $H(v^\lambda(t))$ is a conserved quantity of the equation \eqref{rescaled-cauchy}.

Let $I$ denote the operator
$$ I := I_N^{1-s}$$
using the notation of the previous Section.  We wish to arrange matters so that
\be{hi-init}
H(I u_0^\lambda) \leq \eps^2.
\end{equation}

First consider the kinetic energy:
$$ | \int \frac{1}{8 \pi^2} (I u_0^\lambda)_x^2 | \lesssim \| \partial_x I u_0^\lambda \|_2^2 \lesssim (N^{1-s} \| |\partial_x|^s u_0^\lambda \|_2)^2
= (N^{1-s} \lambda^{-\frac{1}{6}-s} \| |\partial_x|^s u_0 \|_2)^2.$$
Thus if we choose
\be{n-def}
N := C^{-1} \eps^{\frac{2}{1-s}} \lambda^{(\frac{1}{6}+s)/(1-s)}
\end{equation}
for a suitable constant $C$, then we see that the kinetic energy is $\ll \eps^2$.  From this and \eqref{sev-bound} we thus have
\be{h1-bound}
\| I u_0^\lambda \|_{H^1} \leq 10 \eps.
\end{equation}

To deal with the potential energy, we observe the Gagliardo-Nirenberg inequality
\be{gagliardo}
|\int v^5\ dx| \lesssim (\int |v_x|^2\ dx)^{\frac{3}{4}} \|v\|_2^{\frac{7}{2}}.
\end{equation}
Applying this with $v = I u_0^\lambda$ we see from \eqref{sev-bound} that the potential energy is also $\ll \eps^2$.  This gives \eqref{hi-init}.

In the next section we shall prove 
\begin{lemma}\label{iter}
Suppose that one has \eqref{h1-bound}, \eqref{sev-bound}, and
\be{h1-small}
H(I u_0^\lambda) \leq 2\eps^2 
\end{equation}
Then (if $\eps$ is sufficiently small) there exists a unique $Y^s_\lambda$ solution $v^\lambda$ of \eqref{rescaled-cauchy} up to time 1 such that one has the estimates
\be{drift}
|H(I v^\lambda(1)) - H(I v^\lambda(0))| \lesssim \lambda^{0+} N^{-\half}.
\end{equation}
Furthermore, one has 
\be{l2-cons}
\| v^\lambda(1) \|_2 \leq \eps.
\end{equation}
and (if $\lambda$ is sufficiently large depending on $\eps$)
\be{h1-cons}
\| I v^\lambda(1) \|_{H^1} \leq 10 \eps.
\end{equation}
\end{lemma}

This would already give the local well-posedness result for large $H^s$ data.  To obtain the global well-posedness result, we iterate the Lemma (using time translation invariance) and observe from \eqref{hi-init} that one can construct an $H^s$ solution $v^\lambda$ for times
$ \sim \eps^2 \lambda^{0-} N^{\half}$.  Undoing the scaling, we see that we have constructed a solution for time
$$ \sim \eps^2 \lambda^{0-} N^{\half} \lambda^{-3} \gtrsim \eps^C \lambda^{\frac{\frac{1}{6}+s}{2(1-s)} - 3-}$$
by \eqref{n-def}.  If $\frac{5}{6} < s < 1$, then the power of $\lambda$ here is positive, and global well-posedness follows\footnote{If we replace $u^3 u_x$ by $u^k u_x$, then the critical regularity $1/6$ (which appears for instance in \eqref{n-def}) changes to $\half - \frac{1}{k}$.  However, the power $N^{-\half}$ in \eqref{drift} should probably remain unchanged since the exponents in Theorem \ref{main} do not depend on $k$.  If we then repeat the above calculation we are led to the heuristic constraint $s > \frac{13}{14} - \frac{2}{7k}$ for gKdV from this method, although the fact that the $L^2$ norm is not sub-critical for $k \geq 4$ causes technical difficulties in making this heuristic rigorous.} by letting $\lambda \to \infty$.

It remains only to show Lemma \ref{iter}.

\section{Proof of Lemma \ref{iter}}\label{iter-sec}
Let $u_0^\lambda$ obey \eqref{h1-bound}, \eqref{sev-bound}, \eqref{h1-small}.  To construct the solution $v^\lambda$ we repeat the arguments in Section \ref{local-sec}, except we replace the space $Y^s$ by the space $I^{-1} Y^1_\lambda$ defined by
$$ \| v \|_{I^{-1} Y^1_\lambda} := \| I v \|_{Y^1_\lambda}.$$
The space $I^{-1} Y^1_\lambda$ is equivalent to $Y^s_\lambda$, but the constants of equivalence depend on $N$.

As before, it suffices to show that the map
\be{contraction-2}
v^\lambda \mapsto \eta(t) (S^\lambda(t) u^\lambda_0 - \int_0^t S^\lambda(t-t') \P(\P((v^\lambda)^3) v^\lambda_x)(t')\ dt')
\end{equation}
is a contraction on the ball
$$ \{ v^\lambda: \| I v^\lambda \|_{Y^1_\lambda} \lesssim \eps \}.$$
of $I^{-1} Y^1_\lambda$, where 
$$ S^\lambda(t) := \exp(\frac{1}{4\pi^2} t \partial_{xxx})$$
is the free evolution operator on $\R / \lambda \Z$.  Note that $I$ commutes with the $S^\lambda(t)$ as well as the cutoff $\eta(t)$.

By repeating the proof of \eqref{grow} we have
$$ \| \eta(t) S^\lambda(t) u^\lambda_0 \|_{I^{-1} Y^1_\lambda}
= \| \eta(t) S^\lambda(t) I u^\lambda_0 \|_{Y^1_\lambda}
\lesssim \| I u^\lambda_0 \|_{H^1} \lesssim \eps.$$

Also, by repeating the proof of Lemma \ref{invert} we have
$$ \| \eta(t) \int_0^t S^\lambda(t - t') F(t')\ dt' \|_{Y^1_\lambda}
\lesssim \| F \|_{Z^1_\lambda};$$
applying $I$, we obtain
$$ \| \eta(t) \int_0^t S^\lambda(t - t') F(t')\ dt' \|_{I^{-1} Y^1_\lambda}
\lesssim \| IF \|_{Z^1_\lambda}.$$
The contraction then obtains from Corollary \ref{rescaled-I}.

We have thus constructed a function $v^\lambda$ on $\T \times \R$ which satisfies the estimate
\be{iv-bound}
\| I v^\lambda \|_{Y^1_\lambda} \lesssim \eps
\end{equation}
and which solves \eqref{rescaled-cauchy} up to time 1.

Now we show \eqref{drift}.
A computation yields the identity
$$ \partial_t H(v(t)) =
-\int v_t (\frac{1}{4 \pi^2} v_{xx} + \frac{1}{4} v^4)\ dx$$
for arbitrary functions $v(x,t)$.  Since
$$ \int (\frac{1}{4\pi^2} v_{xxx} + v^3 v_x) (\frac{1}{4 \pi^2} v_{xx} + \frac{1}{4} v^4)\ dx
= \int \frac{1}{2} \partial_x (\frac{1}{4 \pi^2} v_{xx} + \frac{1}{4} v^4)^2\ dx = 0
$$
and
$$ \int v_x (\frac{1}{4 \pi^2} v_{xx} + \frac{1}{4} v^4)\ dx = 0$$
and $\P(\P(v^3) v_x) = \P(v^3) v_x$, we see that
$$
\partial_t H(v(t)) =
-\int (v_t + \frac{1}{4\pi^2} v_{xxx} + \P(\P(v^3) v_x))
(\frac{1}{4 \pi^2} v_{xx} + \frac{1}{4} v^4)\ dx.
$$
Now applying this with $v := I v^\lambda$.  From \eqref{rescaled-cauchy}, and the fact that $I$ commutes with derivatives and $\P$, we have
$$ (I v^\lambda)_t + \frac{1}{4\pi^2} (I v^\lambda)_{xxx} + 
\P I(\P((v^\lambda)^3) v^\lambda_x)) = 0.$$
Inserting this into the previous, we can express $\partial_t H(I v^\lambda(t))$ as a commutator:
$$
\partial_t H(I v^\lambda(t)) =
\int \P
\left(
I(\P(v^\lambda v^\lambda v^\lambda) v^\lambda_x) - 
\P(I v^\lambda I v^\lambda I v^\lambda) I v^\lambda_x
\right)
 F
\ dx
$$
where
$$ F := \frac{1}{4 \pi^2} Iv^\lambda_{xx} + \frac{1}{4} (Iv^\lambda)^4.$$

From the fundamental theorem of Calculus, it thus suffices to show that
$$
|\int\int \chi_{[0,1]}(t) \P
\left(
I(\P(v^\lambda v^\lambda v^\lambda) v^\lambda_x) - 
\P(I v^\lambda I v^\lambda I v^\lambda) Iv^\lambda_x
\right)
F\ dx dt| \lesssim \lambda^{0+} N^{-\half+}.$$
By Lemma \ref{duality} it suffices to show the estimates
\be{f-bound}
\| F \|_{Y^{-1}_\lambda} \lesssim \lambda^{0+}
\end{equation}
and
\be{commutator}
\| \P\left(
I(\P(v^\lambda v^\lambda v^\lambda) v^\lambda_x) - 
\P(I v^\lambda I v^\lambda I v^\lambda) I v^\lambda_x
\right)
\|_{Z^1_\lambda} \lesssim \lambda^{0+} N^{-\half}.
\end{equation}

We first show \eqref{f-bound}.  The contribution of the main term $I v^\lambda_{xx}$ is acceptable from \eqref{iv-bound}.  To control the lower order term $(Iv^\lambda)^4$ it suffices by \eqref{iv-bound} to show that
$$
\| v^4 \|_{Y^{-1}_\lambda} \lesssim 
\lambda^{0+} \| v \|_{Y^1_\lambda}^4$$
for all functions $v$.  This estimate will be obtained with plenty of room to spare, given that we are $\frac{7}{6}$ derivatives above scaling.

From Proposition \ref{main-lambda} with $k=4$ we have
$$
\| v^4 \|_{X^{-\half, \half}_\lambda} \lesssim 
\lambda^{0+} \| v \|_{Y^1_\lambda}^4,$$
so it suffices by \eqref{y-def} (and discarding a derivative) to show that
$$
\| \widehat{v^4}(\xi,\tau) \|_{L^2_\xi L^1_\tau } \lesssim 
\| \langle \xi \rangle \hat v \|_{L^2_\xi L^1_\tau}^4,$$
which is of course equivalent to
$$
\| (\hat v * \hat v * \hat v * \hat v)(\xi,\tau) \|_{L^2_\xi L^1_\tau } \lesssim 
\| \langle \xi \rangle \hat v \|_{L^2_\xi L^1_\tau}^4.$$
Since convolutions of $L^1_\tau$ functions stay in $L^1_\tau$, it suffices to show the spatial estimate
$$ \| (F_1 * F_2 * F_3 * F_4)(\xi) \|_{L^2_\xi} \lesssim 
\prod_{i=1}^4 \| \langle \xi \rangle F_i \|_{L^2_\xi},$$
which by Plancherel is equivalent to
$$ \| f_1 f_2 f_3 f_4 \|_{L^2_x} \lesssim \prod_{i=1}^4 \|f_i\|_{H^1_x}.$$
But this follows from H\"older and the Sobolev embedding $H^1_x \subseteq L^8_x$.

It remains to show \eqref{commutator}.  By \eqref{iv-bound} it suffices to show
\be{com}
\| \P(I(\P(v_1 v_2 v_3) (v_4)_x)) - \P((I v_1) (I v_2) (Iv_3)) (Iv_4)_x \|_{Z^1_\lambda} \lesssim \lambda^{0+} N^{-\half}
\prod_{i=1}^4 \| I v_i \|_{Y^1_\lambda}
\end{equation}
for all $v_1, \ldots, v_4$.  This will be accomplished by Corollary \ref{rescaled} (at the endpoint $s=\half$) and a variant of the arguments in Lemma \ref{interp-lemma}.  In order to obtain the crucial factor of $N^{-\half+}$ we must exploit the cancellation between the two terms on the left-hand side of \eqref{com}.

We turn to the details.  Without loss of generality we may assume that $v_1$, $v_2$, $v_3$, and $(v_4)_x$ have non-negative Fourier transforms.  We divide into three cases.

{\bf Case 1:  (Low-low interactions)} $v_1$, $v_2$, $v_3$, and $(v_4)_x$ all have Fourier support in the region $|\xi| \leq \frac{N}{5}$.

In this case all the $I$ operators act like the identity, and the left-hand side of \eqref{com} vanishes.

{\bf Case 2:  (High-high interactions)} At least two of $v_1, v_2, v_3, (v_4)_x$ have Fourier support in the region $|\xi| \geq \frac{N}{100}$.

In this case we will not exploit any cancellation in \eqref{com}.  From the observation that
$$ m^{1-s}( \frac{\xi_1}{N}) m^{1-s}( \frac{\xi_2}{N}) m^{1-s}( \frac{\xi_3}{N}) m^{1-s}( \frac{\xi_4}{N})
\lesssim m^{1-s}( (\xi_1 + \xi_2 + \xi_3 + \xi_4) / N)$$
for all $\xi_1, \ldots, \xi_4 \in \R$, we see that it suffices to show that
$$
\| I\P(\P(v_1 v_2 v_3) (v_4)_x) \|_{Z^1_\lambda} \lesssim 
\lambda^{0+} N^{-\half+} \prod_{i=1}^4 \| I v_i \|_{Y^1_\lambda}.$$
We now use the estimate
$$ \| I v \|_{Z^1_\lambda} \lesssim N^{1-s} \| \langle \nabla \rangle^{s-1} v \|_{Z^1_\lambda}
\sim N^{1-s} \| \langle \nabla \rangle^{s-\half} v \|_{Z^{\half}_\lambda},$$
where $\langle \nabla \rangle$ is the Fourier multiplier with symbol $\langle \xi \rangle$, to estimate the left-hand side by
$$ N^{1-s}
\| \langle \nabla \rangle^{s-\half} \P((\P(v_1 v_2 v_3) (v_4)_x)) \|_{Z^{\half}_\lambda}.$$
Applying the fractional Leibniz rule (using the positivity of the Fourier transforms of $v_1, v_2, v_3, (v_4)_x$), we may distribute the differentiation operator $\langle \nabla \rangle^{s-\half}$ to one of the functions, say $(v_4)_x$.  (The other cases are similar and will be left to the reader).  We thus estimate the previous by
$$ N^{1-s}
\| \P((\P(v_1 v_2 v_3) (\langle \nabla \rangle^{s-\half} v_4)_x)) \|_{Z^{\half}_\lambda}.$$
Applying Corollary \ref{rescaled} with $s=\half$, we can estimate this by
$$ \lambda^{0+} N^{1-s}
\| v_1 \|_{Y^{\half}_\lambda} \| v_2 \|_{Y^{\half}_\lambda}
\| v_3 \|_{Y^{\half}_\lambda} \| \langle \nabla \rangle^{s-\half} v_4 \|_{Y^{\half}_\lambda}.$$
Note that 
$$\| \langle \nabla \rangle^{s-\half} v_4 \|_{Y^{\half}_\lambda} \lesssim 
N^{s-1} \| I v_4 \|_{Y^1_\lambda}.$$
Also, at least one of the functions $v_1$, $v_2$, $v_3$ has Fourier support in the region $|\xi| \geq N/100$, so for this function we have
$$ \| v_i \|_{Y^{\half}_\lambda} \lesssim N^{-1/2} \| I v_i \|_{Y^1_\lambda}.$$
For the other two functions we just use the crude bound
\be{crude}
\| v_i \|_{Y^{\half}_\lambda} \lesssim \| I v_i \|_{Y^1_\lambda}.
\end{equation}

Combining all these estimates we obtain the result.

The only remaining case to consider is 

{\bf Case 3: (Low-high interactions)} One of $v_1, v_2, v_3, (v_4)_x$ has Fourier support in the region $|\xi| \geq \frac{N}{5}$, and the other three have Fourier support in the region $|\xi| \leq \frac{N}{100}$.

Let us suppose that $(v_4)_x$ is the function with Fourier support in the region $|\xi| \geq \frac{N}{5}$; the other cases are similar.  The idea will be to exploit the cancellation in \eqref{com} to transfer one derivative from the high-frequency function $v_4$ to the low frequency functions $v_1, v_2, v_3$.

In this case the operator $I$ is the identity on $v_1, v_2, v_3$, so we may write the left-hand side of \eqref{com} as
\be{com-2}
\| \P(I(\P(v_1 v_2 v_3) (v_4)_x)) - \P(v_1 v_2 v_3) (Iv_4)_x) \|_{Z^1_\lambda}.
\end{equation}
From the mean-value theorem we observe that
$$ | m^{1-s}( \xi' + \xi ) - m^{1-s}( \xi ) | \lesssim \frac{ |\xi'| }{ |\xi| } m^{1-s} (\xi )$$
if $|\xi| \geq \frac{1}{5}$ and $|\xi'| \leq \frac{4}{100}$.  Thus we have
\begin{align*}
 | m^{1-s}( (\xi_1 + \xi_2 + \xi_3 + \xi_4)/N ) - m^{1-s}( \frac{\xi_4}{N} ) | &\lesssim \frac{ |\xi_1 + \xi_2 + \xi_3| }{ |\xi_4| } m^{1-s} (
\frac{\xi_4}{N} ) \\
&\lesssim \frac{ |\xi_1 + \xi_2 + \xi_3|^{1/2} }{ N^{1/2} } \frac{N^{1/2}}{|\xi_4|^{1/2}} m^{1-s} (
\frac{\xi_4}{N} ) \\
&\lesssim N^{-1/2} (|\xi_1|^{1/2} + |\xi_2|^{1/2} + |\xi_3|^{1/2}) m^{3/2-s} (
\frac{\xi_4}{N} ) 
\end{align*}
if $|\xi_4| \geq N/5$ and $|\xi_1|, |\xi_2|, |\xi_3| \leq \frac{N}{100}$. We may thus estimate \eqref{com-2} by
$$
N^{-1/2} \| \P(\P((|\nabla|^{1/2} v_1) v_2 v_3) (I_N^{3/2-s} v_4)_x) \|_{Z^1_\lambda}
\sim
N^{-1/2} \| \langle \nabla \rangle^{\half} \P(\P((|\nabla|^{1/2} v_1) v_2 v_3) (I_N^{3/2-s} v_4)_x) \|_{Z^{\half}_\lambda}
$$
together with three other similar terms.
Since $v_4$ has higher frequency than the other three functions, we may distribute the derivative $\langle \nabla \rangle^{\half}$ onto $v_4$.  By Corollary \ref{rescaled} with $s=\half$, we can estimate the previous by
$$ \lambda^{0+} N^{-1/2}
\| |\nabla|^{1/2} v_1 \|_{Y^{\half}_\lambda} \| v_2 \|_{Y^{\half}_\lambda}
\| v_3 \|_{Y^{\half}_\lambda} \| \langle \nabla \rangle^{\half} I_N^{3/2-s} v_4 \|_{Y^{\half}_\lambda}.$$
From the frequency support of $v_4$ we see that
$$ \| \langle \nabla \rangle^{\half} I_N^{3/2-s} v_4 \|_{Y^{\half}_\lambda} \lesssim 
\| I v_4 \|_{Y^1_\lambda},$$
and the claim follows from this, noting from the frequency support of $v_1,v_2,v_3$ that
$$\| |\nabla|^{1/2} v_i \|_{Y^{\half}_\lambda}, \| v_i \|_{Y^{\half}_\lambda}  \lesssim \|v_i\|_{Y^1_\lambda}.$$
This completes the proof of \eqref{com} in all cases, and \eqref{drift} follows.

The estimate \eqref{l2-cons} follows from \eqref{sev-bound} and the observation that the $L^2$ norm is conserved by the flow \eqref{rescaled-cauchy}.  Now we show \eqref{h1-cons}.  By \eqref{l2-cons} it suffices to show that
\be{couple}
\| \partial_x I v^\lambda(1) \|_{L^2_x} \leq 5 \eps.
\end{equation}

From \eqref{iv-bound}, \eqref{sup-est} we have
$ \| Iv^\lambda(1) \|_{H^1} \lesssim \eps;$
from \eqref{gagliardo} we thus have
$ |\int (I v^\lambda(0))^5\ dx| \lesssim \eps^5.$
The claim \eqref{couple} then follows from \eqref{h1-small}, \eqref{drift} if $\eps$ is sufficiently small and $\lambda$ (and hence $N$) is sufficiently large.
\endprf

\end{document}